\documentclass[10pt, a4paper, twoside]{amsart}

\usepackage[latin1]{inputenc}

\usepackage{stmaryrd}

\usepackage{amssymb}
\usepackage{array}

\usepackage{amsmath}
\usepackage{amsthm}
\usepackage{amscd}
\usepackage{delarray}

\usepackage{mathrsfs}

    \makeatletter
    \newcommand{\cdotsfor}[1]{%
      \ifx[#1\@xp\scdots@for\else\cdots@for\@ne{#1}\fi}
    \newmuskip\dotsspace@
    \def\scdots@for#1]{\cdots@for{#1}}
    \def\cdots@for#1#2{\multicolumn{#2}c%
      {\m@th\dotsspace@1.5mu\mkern-#1\dotsspace@
       \xleaders\hbox{$\m@th\mkern#1\dotsspace@\cdot\mkern#1\dotsspace@$}%
               \hfill
       \mkern-#1\dotsspace@}%
       }
    \makeatother

    \makeatletter  
    \def\adots{\mathinner{\mkern2mu\raise\p@\hbox{.}
    \mkern2mu\raise4\p@\hbox{.}\mkern1mu
    \raise7\p@\vbox{\kern7\p@\hbox{.}}\mkern1mu}}
    \makeatother

\newcommand{\bbA}{\mathbb{A}}  \newcommand{\bbC}{\mathbb{C}}  \newcommand{\bbF}{\mathbb{F}}  
\newcommand{\bbG}{\mathbb{G}}  \newcommand{\bbN}{\mathbb{N}}  \newcommand{\bbP}{\mathbb P}
    \newcommand{\bbZ}{\mathbb{Z}}

    \newcommand{\C}{\bbC}
    \newcommand{\N}{\bbN}
    
    \newcommand{\Z}{\bbZ}
    \renewcommand{\P}{\bbP}

  \newcommand{\cE}{\mathcal E}  \newcommand{\cO}{\mathcal O}

   \newcommand{\E}{\cE}

\newcommand{\ssA}{\mathsf{A}}

\newcommand{\ssD}{\mathsf{D}}

\newcommand{\bF}{\mathbf{F}}

\newcommand{\hz}[1]{H^0({#1})}

\newcommand{\donne}{\Longrightarrow}

\newcommand{\vers}{\longrightarrow}

\newcommand{\Spc}[1]{\langle {#1} \rangle}
\newcommand{\pe}[2]{\stackrel{#1}{#2}}

\newcommand{\tens}{\otimes}

\renewcommand{\leq}{\leqslant}  
\renewcommand{\geq}{\geqslant}  

\newcommand{\blt}{\bullet}  

\newcommand{\cst}{\C^{*}}

\theoremstyle{plain} 
\newtheorem{thm}{Theorem}[section]
\newtheorem*{thm*}{Theorem}
\newtheorem{prop}[thm]{Proposition}
\newtheorem{lema}[thm]{Lemma}
\newtheorem{cor}[thm]{Corollary}

\theoremstyle{definition} 
\newtheorem{defin}[thm]{Definition}
\newtheorem{numar}[thm]{}  

\newtheorem{notation}[thm]{Notation}

\theoremstyle{remark}
\newtheorem{remarca}[thm]{Remark}
\newtheorem*{ack}{Acknowledgements}
\newtheorem*{overview}{Overview}

\DeclareMathOperator{\codim}{codim}
  
\DeclareMathOperator{\Ker}{Ker}

\DeclareMathOperator{\Aut}{Aut}  
\DeclareMathOperator{\card}{card}

\DeclareMathOperator{\diag}{diag}

\DeclareMathOperator{\slin}{\mathfrak{sl}}  
\DeclareMathOperator{\glin}{\mathfrak{gl}}

\DeclareMathOperator{\Span}{Span}

\DeclareMathOperator{\rank}{rank}

\newcommand{\GL}{\mathrm{GL}}  
\newcommand{\SL}{\mathrm{SL}}  
  
\newcommand{\PSL}{\mathrm{PSL}}  
\newcommand{\Sp}{\mathrm{Sp}}  
\newcommand{\SSp}{\mathrm{SSp}}  
\newcommand{\PSp}{\mathrm{PSp}}  
\newcommand{\liesp}{\mathfrak{sp}}  
\newcommand{\liespd}{\liesp_{2n}}  
\newcommand{\liespdd}{\liesp_{2n+1}}  
\newcommand{\liespddd}{\liesp_{2n+2}}  
\newcommand{\lietspdd}{\widetilde{\liesp}_{2n+1}}

\newcommand{\spd}{\Sp_{2n}}  
\newcommand{\spdd}{\Sp_{2n+1}} 
\newcommand{\spddd}{\Sp_{2n+2}}
  
\newcommand{\sspdd}{\SSp_{2n+1}}

\newcommand{\pspd}{\mathrm{PSp}_{2n}}   
\newcommand{\pspdd}{\mathrm{PSp}_{2n+1}}

\newcommand{\csp}{\C^*\times\Sp_{2n}} 
  
\newcommand{\gld}{\GL(2n,\C)}  
\newcommand{\gldd}{\GL(2n+1,\C)}  
  
\newcommand{\tldo}{\tilde\omega }   
\newcommand{\tlde}{\tilde E}

\newcommand{\transp}[1]{{\,}^t\!#1}   

\newcommand{\gokd}{G_{\omega }(k,2n)}  
\newcommand{\gokdd}{G_{\omega }(k,2n+1)}  
\newcommand{\gokddd}{G_{\omega }(k,2n+2)}  
\newcommand{\gokE}{G_{\omega }(k,E)}  
\newcommand{\goke}{\gokE}  
\newcommand{\go}{G_{\omega}}

\newcommand{\gke}{G(k,E)}

\newcommand{\sldd}{\SL_{2n+1}}

\newcommand{\indiceadm}[2]{\mathbf{I}^{\omega }_{#1,#2}   }  
\newcommand{\iod}{\indiceadm{k}{2n}}
\newcommand{\iodd}{\indiceadm{k}{2n+1}}
\newcommand{\ioddd}{\indiceadm{k}{2n+2}}
\newcommand{\indice}[2]{\mathbf{I}_{#1,#2}}  
\newcommand{\idn}{\indice{k}{2n}}   
\newcommand{\idd}{\indice{k}{2n+1}}

     \newcommand{\resdst}[2]{{{#1}\raisebox{-.23ex}{$|$}}{}^{\phantom{,}}_{#2}}
     \newcommand{\ressst}[2]{{{#1}\raisebox{-.23ex}{$\scriptstyle|$}}{}^{\phantom{,}}_{#2}}
     \newcommand{\resssst}[2]{{{#1}\raisebox{-.25ex}{$\scriptscriptstyle|$}}{}^{\phantom{,}}_{#2}}
\newcommand{\restst}[2]{{{#1}\raisebox{-.20ex}{$|$}}{}_{#2}}

\newcommand{\res}[2]{\mathchoice{\resdst{#1}{#2}}{\restst{#1}{#2}}{\ressst{#1}{#2}}{\resssst{#1}{#2}}}  

\newcommand{\ldt}{\Lambda ^2T^*}

\newcommand{\mx}{\mathbf{max}}  
\renewcommand{\b}[1]{\overline{#1}} 

\newcommand{\e}{\varepsilon }

\newcommand{\Id}{\mathrm{Id}}

\newcommand{\liee}{\mathfrak{e}}

\newcommand{\tast}{(T^a)^*}  
\newcommand{\tbst}{(T^b)^*}

\newcommand{\fo}{\bbF_{\omega}}  
\newcommand{\fod}{\bbF_{\omega} (2n)}  
\newcommand{\fodd}{\bbF_{\omega} (2n+1)}  
\newcommand{\foddd}{\bbF_{\omega} (2n+2)}

\newcommand{\foE}{\bbF_\omega (E)}  
  
\newcommand{\foe}{\foE}

\newcommand{\cd}{\C^{2n}}
\newcommand{\cdd}{\C^{2n+1}}
\newcommand{\cddd}{\C^{2n+2}}

\newcommand{\liep}{\mathfrak{p}}
\newcommand{\liet}{\mathfrak{t}}
\newcommand{\lieu}{\mathfrak{u}}
\newcommand{\lieb}{\mathfrak{b}}

\newcommand{\liesl}{\mathfrak{sl}}

\newcommand{\liesldd}{\liesl_{2n+1}}

\newcommand{\liegl}{\mathfrak{gl}}
\newcommand{\liegld}{\liegl_{2n}}

\newcommand{\wspd}{W(\spd)}  
\newcommand{\wspdd}{W(\spdd)}  
\newcommand{\wspddd}{W(\spddd)}

\newcommand{\tSp}{\widetilde{\mathrm{Sp}}}  

\newcommand{\tspdd}{\tSp_{2n+1}}

\newcommand{\sspc}[1]{S\raisebox{-.35ex}{$\scriptstyle\langle #1 \rangle$}} 
\newcommand{\ssl}{\sspc{\lambda }}

\newcommand{\tddd}{T_{2n+2}}  

\newcommand{\lietd}{\liet_{2n}}  
  
\newcommand{\lietddd}{\liet_{2n+2}}

\newcommand{\bddd}{B_{2n+2}}  

\newcommand{\liebd}{\lieb_{2n}}  
\newcommand{\liebdd}{\lieb_{2n+1}}  
\newcommand{\liebddd}{\lieb_{2n+2}}

\newcommand{\uplsddd}{U^+_{2n+2}}

\newcommand{\lieuplsd}{\lieu^+_{2n}}  
\newcommand{\lieuplsdd}{\lieu^+_{2n+1}}  
\newcommand{\lieuplsddd}{\lieu^+_{2n+2}}  

\newcommand{\lieumnsd}{\lieu^-_{2n}}  
\newcommand{\lieumnsdd}{\lieu^-_{2n+1}}  
\newcommand{\lieumnsddd}{\lieu^-_{2n+2}}

\usepackage[arrow, curve, matrix, rotate]{xy}
\xyoption{dvips}
\xyoption{ps}  

\newcommand{\bari}{\bar{\imath}}  
\newcommand{\barj}{\bar{\jmath}}  

\newcommand{\tE}{\tilde{E}}

\newcommand{\fto}{\bbF_{\tldo}}
\newcommand{\wz}{\bar{1}\bar{2}\dots\bar{n}0}

\begin{document}

\title{Odd symplectic flag manifolds} 
\author{Ion Alexandru Mihai} 

\begin{abstract}
We define the \emph{odd symplectic} grassmannians and flag manifolds,
which are smooth projective varieties 
equipped with an action of the odd symplectic group and
generalizing the usual symplectic grassmannians and flag manifolds.
Contrary to the latter, which are the flag manifolds of the
symplectic group, the varieties we introduce are not homogeneous. 
We argue nevertheless that in 
many 
respects the odd symplectic grassmannians and flag manifolds
behave like homogeneous varieties; in support of this claim, 
we compute the automorphism group
of the odd symplectic grassmannians, and we prove a Borel-Weil type
theorem for the odd symplectic group.
\end{abstract}

\maketitle

\section{Introduction}
\label{sec:intro}

In \cite{Proctor1} Proctor introduces the \emph{odd symplectic
  group}, a generalization of the symplectic group on an odd
dimensional space defined as the group of linear transformations
preserving a generic skew-form.
His initial motivation is a series of combinatorial identities, but
Proctor goes on to study a certain class of representations of this
group (which is not reductive) 
and eventually proves a character formula very similar to Weyl's
formula for the simple Lie groups. 
In this way the odd symplectic group 
presents some similarities with the simple Lie groups and
appears to fit nicely in the framework of the classical groups, filling 
the ``gap'' in the series $ \{\spd \}_{n}$.

Here, we take a different look at this situation, from a more geometric
perspective. Recall that the flag manifolds of the symplectic group,
$\spd /P$, with $P$ a parabolic subgroup in $\spd$, identify with the
varieties of flags of isotropic subspaces of $\cd$.
We generalize these varieties to the odd
symplectic situation, in the most 
straightforward way : let $\omega $ be a generic skew-form on $\cdd$
and define the \emph{odd symplectic flag manifolds} 
to be the varieties of flags of subspaces of $\cdd$ 
isotropic with respect to $\omega $. 
These are projective varieties 
and are
equipped with natural actions of the odd symplectic group $\spdd$  
which preserves the
skew-form $\omega $. Unlike in the symplectic setting however, these actions
are no longer transitive, as the non trivial kernel of the skew-form
is preserved by $\spdd$ and therefore isotropic flags
having different incidence relations with this kernel cannot be in the
same orbit.

The aim of this paper is to present some evidence that 
the odd symplectic flag manifolds, although not homogeneous, 
fill the role of flag manifolds for the odd symplectic group.
Part of this evidence, for example, is constituted by the theorems
below describing respectively the automorphism group of the odd
symplectic grassmannians and a Borel-Weil theorem for the odd
symplectic group, which show 
that sometimes the symplectic flag
manifolds and their odd symplectic counterparts behave like a ``series''.

\medskip

The reason we do not simply consider the
homogeneous spaces $\spdd/P$, for $P \subset \spdd$ a parabolic
  subgroup, as ``flag manifolds'' for $\spdd$, is that, as it will become clear below, these
  coincide with the flag 
  manifolds of the symplectic group $\spd$,
and therefore do not constitute representative examples 
for the odd symplectic situation.

The odd symplectic flag manifolds are not homogeneous but it turns
out they are quasi-homogeneous. The odd symplectic group acts with
finitely many orbits, described by the incidence relations with the
kernel
mentioned above. 
They are also smooth, which follows from the fact that, as subvarieties of
the (type $\ssA$) flag manifolds defined by forgetting the isotropy
conditions, they appear as the zero locus of a generic section of a
vector bundle.

A natural question is whether the odd symplectic flag manifolds admit
a cellular decomposition similar to the decomposition into Schubert
cells of the usual symplectic flag manifolds. In the symplectic case,
and 
in general for all the flag manifolds of the classical groups, the Schubert
cells can be defined in several equivalent ways, eg 
they can be described
by incidence relations with respect to a fixed flag and they coincide
with the orbits of a Borel subgroup.
It turns out that these two 
recipes can be used in the odd
symplectic setting as well to define cell decompositions
of the odd symplectic
flag manifolds.

Actually, the odd symplectic flag manifolds
themselves can be identified with certain Schubert subvarieties in
symplectic flag manifolds 
(and then the cell decompositions above coincide
with their cell decomposition as Schubert varieties).
This goes as follows: the generic skew-form $\omega $ on $\cdd$ can be
extended to a symplectic form $\tldo$ on $\cddd$, so that any odd
symplectic flag manifold associated to $\omega $ is identified with
the Schubert subvariety of the corresponding symplectic flag manifold
associated to $\tldo$ given by those flags which are contained in the
hyperplane $\cdd$.
The parabolic subgroup of $\spddd$ which preserves the hyperplane
$\cdd$ acts therefore on the odd symplectic flag manifolds, via the
morphism of restriction to $\cdd$ which is surjective with image
$\spdd$.

The subgroup $\tspdd$ of this parabolic subgroup which fixes an equation of the
hyperplane $\cdd$ (or, equivalently, a vector $e$ generating the
kernel of $\omega $) has been considered before by Gelfand and
Zelevinski in \cite{GelfZelev} as a  variant odd symplectic group, in
connection to the problem of constructing representation models for
the classical groups. 
We will call it the \emph{intermediate odd symplectic group} to
distinguish it from the odd symplectic group considered above. 
This group has been considered later by Shtepin in \cite{Shtepin} 
where he constructs a series of $\tspdd$-modules as a
means to separate multiple components 
when restricting simple $\spddd$-modules to $\spd$.
Its Lie algebra $\lietspdd$ coincides with the \emph{intermediate} Lie
algebra of the symplectic Lie algebra $\liespddd$, a general
construction which can be associated to any simple Lie algebra (which
is used for example in \cite{Sextonions} to construct the Lie algebra
$\liee_{7\frac12}$ sitting between the exceptional Lie algebras
$\liee_7$ and $\liee_8$).

\medskip

In this paper we will focus only on the extremal types of odd
symplectic flag manifolds, namely the \emph{odd symplectic grassmannians}
$$
\gokdd =  \{ V \mid V\subset \cdd, \dim V=k, V\text{ isotropic}/\omega   \} 
$$ 
and the variety of maximal flags of isotropic
subspaces $\fodd$, which we
simply call the \emph{odd symplectic flag manifold}.
Since the maximal isotropic subspaces in $\cdd$ are
those of dimension $n+1$, the variety $\fodd$ is the variety of flags
of the form
$V_{\blt}=(V_1\subset V_2\subset\dots\subset V_{n+1})$
with each $V_i \subset \cdd$ isotropic of dimension $i$.  
These varieties are the analogues in the odd symplectic setting of the
\emph{symplectic grassmannians} $\gokd$ and the 
\emph{odd symplectic flag manifold} $\fod$, 
which are, respectively, the minimal and the maximal flag varieties of
the symplectic group $\spd$.

We compute the automorphism group of the odd symplectic grassmannians
$\gokdd$
and find out that for $2\leq k \leq n$ it equals $\pspdd$, the quotient of $\spdd$ by its
center $ \{\pm 1 \}$ (for $k=n+1$ the odd symplectic grassmannian 
$\go(n+1,2n+1)$ is isomorphic to the symplectic grassmannian
$\go(n,2n)$ and therefore its automorphism group is $\pspd$).
In this respect the odd symplectic grassmannians behave like 
homogeneous flag manifolds, 
since, as a general rule, the connected
automorphism group $\Aut^{\circ} (G/P)$ of a flag manifold is the
adjoint group $\mathrm PG$ (there are some exceptions to this rule,
see \cite[\S3.3, Theorem 2]{Akh}) .
This is also suggestive of the behavior of the symplectic and odd
symplectic grassmannians as a series, as we can state this result in
the following uniform manner:

\begin{thm*} 
For integers $N$ and $k$ such that $2 \leq k \leq [N/2]$, the automorphism group of the variety 
$\go(k,N)$ is $\PSp_N=\Sp_N / \{\pm1\}$. 
\end{thm*}

Finally, this also shows a 
close connection between the group $\spdd$ and 
the odd symplectic grassmannians $\gokdd$ since, up to the center
$ \{\pm 1 \}$, we can recover $\spdd$
from the
geometry of
$\gokdd$ (which is not the case for the homogeneous spaces $\spdd/P$,
$P$ a parabolic subgroup of $\spdd$).

Another natural question is whether there is an analog of the
Borel-Weil theorem for the odd symplectic group. 
The usual Borel-Weil theorem connects  the representation theory of a
simple complex 
Lie
group $G$ with the geometry of the flag variety
$G/B$ by explicitly 
identifying the simple $G$-modules with the spaces of global sections
of the line bundles on $G/B$.
There are two points that need to be addressed in order to generalize
this to the odd symplectic setting. 
First, the odd symplectic group is not reductive and therefore the
simple modules do not play the same role in its representation
theory
as they do in the symplectic case.
We need to replace
them by another class of preferred representations, 
and it is natural to consider 
the class of $\spdd$-modules introduced by Proctor.
These are defined by porting to the odd symplectic setting 
the construction of Weyl
of the simple modules of the symplectic
group. Specifically, for $\lambda $ a partition with at most $n+1$
parts define the $\spdd$-module $\sspc{\lambda }\cdd$ as the
intersection of the Schur power $S_\lambda \cdd$ with the kernels of
all the possible contractions with the odd symplectic form $\omega $
(the ``trace free'' part of $S_\lambda \cdd$).
We may also
consider the representations of the intermediate odd
symplectic group $\tspdd$ which were introduced by Shtepin. 
We show that actually these are isomorphic to the representations
defined by Proctor, via the natural morphism $\tspdd \to \spdd$.
The second point that needs attention is that
the odd symplectic flag manifold $\fodd$ is not homogeneous,
so we no longer have a correspondence between line bundles and
characters as in the symplectic case. 
It is natural then to use the fact that $\fodd$ identifies with a Schubert
subvariety in the symplectic flag manifold $\foddd$ and consider 
those line bundles on $\fodd$ which come from $\foddd$. 
We write them
in terms of the tautological bundles, and eventually obtain:

 \begin{thm*} 
For $N$ an integer let $n=[(N-1)/2]$. 
Let $\lambda =(\lambda _0 \geq \lambda _1 \geq \cdots \geq \lambda _n)$
be a partition with at most $n+1$ parts. Denote by $L_\lambda $ the
line bundle on  $\fo(\C^N)$ given by
$$
L_\lambda  = 
{T_1^*}^{\tens \lambda _0}    \tens
    {(T_2/T_1)^*}^{\tens \lambda _1}    \tens\cdots\tens 
    {(T_{n+1}/T_n)^*}^{\tens \lambda _n}
$$
 where $T_i$ is the rank $i$ tautological vector bundle on
 $\fo(\C^N)$.
Then, as $\Sp_N$-modules, we have
$$
H^0(\fo(\C^N), L_\lambda) \simeq (\sspc{\lambda }\C^N)^*.
$$
 \end{thm*}

Here, when $N=2n+2$ we get the usual Borel-Weil theorem for the
symplectic group $\Sp_{2n+2}$. 
So again we notice the behavior of the symplectic and odd symplectic
flag manifolds as a series.

\medskip

In a forthcoming paper (\cite{eqcoh}) we study the equivariant
cohomology of the odd symplectic flag manifold $\fodd$ correspondig to
the action of a maximal torus of $\spdd$. In particular, we compute
the singular cohomology algebra $H^*(\fodd, \C)$ which turns out to be
isomorphic 
to the singular cohomology algebra of the
flag manifold $G/B$ for $G$ of type $\ssD_{n+1}$. 

\medskip

\begin{overview} 
This paper is organized as follows. In the second section we gather
some basic facts and fix notation concerning the symplectic groups and
their flag manifolds. In the third section we introduce  Proctor's
odd symplectic group and list some of its properties. We also recall
here the intermediate odd symplectic group of Gelfand and Zelevinski
and its relation with the odd symplectic group.
In the fourth section we define the odd symplectic flag manifolds. We
describe the $\spdd$-orbits and the Schubert decompositions in the
case of the odd symplectic grassmannians and the odd symplectic flag
manifold. The fifth section deals with the computation of the
automorphism group of the odd symplectic grassmannians. We conclude,
in the sixth section, with the Borel-Weil theorem for the odd
symplectic group.
\end{overview}

\begin{ack}
I warmly thank Laurent Manivel for his invaluable help and support in
writing this article. I also thank Patrick Polo for his very useful
comments concerning section $6$. This paper has been written while I
was a post-doc fellow at the Weizmann Institute of Science under the
Marie Curie RTN \emph{Liegrits}. I wish to thank my hosts for their
hospitality and the perfect working conditions they provided.
\end{ack}

\tableofcontents

\section{Preliminaries}

We start by reviewing some basic facts about the symplectic group and
its flag manifolds.

\begin{numar} 
Let $F$ complex vector space of dimension $2n$ and $\omega \in \Lambda ^2F^*$ 
a \emph{symplectic} form on $F$, ie a skew-symmetric,
non-degenerate form.
A subspace $V\subset F$ is \emph{isotropic} if 
$\omega (x,y)=0$ for all $ x,y\in V$.
A \emph{symplectic basis} is a basis $\{e_1, \dots, e_{2n}\}$  of $F$ 
such that
$$
\omega (e_i, e_j)=\delta _{i,2n+1-j}
$$
for all $1\leq i,j \leq 2n$.
Any symplectic form admits a symplectic basis.
A standard notation convention, which we will also use in this paper,
is to denote 
$\bar{\imath}=2n+1-i$ for 
$i\in\{1,\dots,2n\}$.
A symplectic basis is therefore a basis
$\{e_1, \dots, e_n, e_{\bar{n}},\dots, e_{\bar{1}} \}$
such that 
$\omega (e_i,e_{\bar{\jmath} })=\delta _{ij}$ 
for all $1\leq i,j \leq 2n$. 
In a symplectic basis, the matrix of the form $\omega $ is
$$
J=
\begin{pmatrix} 0 & A \\  -A & 0
\end{pmatrix}
$$
where $A$ is the $n \times n$ matrix whose anti-diagonal entries are
all equal to 1 and all the other entries are 0.
\end{numar}

\begin{numar} 
The \emph{symplectic group} $\spd$
is the stabilizer of $\omega $ in $\GL(F)$ for the standard action of
$\GL(F)$ in $\Lambda ^2F^*$.
It is a simple algebraic group.
The choice of a symplectic basis identifies $\spd$ with the matrix group
$$
\{g\in \gld \mid  \transp{g} J g = J \}.
$$
From now on, we fix a symplectic basis 
$\{e_1, \dots, e_{2n}\}$ 
and identify $\spd$ with this matrix group.

The subgroup of $\spd$ of upper triangular matrices in
the 
basis
$\{ e_1, \dots, e_{2n} \}$ is a Borel subgroup.
Also, the subgroup $T_{2n}$  of $\spd$ of diagonal matrices in the
basis $\{ e_1, \dots, e_{2n}  \}$ is a maximal torus.
We have
$$T_{2n} =\{\diag(t_1, \dots, t_n, t_n^{-1}, \dots, t_1^{-1}) \mid t_1, \dots, t_n \in \cst\},$$
where $\diag(x_1,\dots, x_{2n})$ denotes 
the diagonal matrix whose diagonal entries are
$x_1,\dots, x_{2n}$.
Denote by $\e_1, \dots, \e_{2n}$ the characters 
$\e_i: \diag(x_1,\dots, x_{2n}) \mapsto x_i$
of $T_{2n}$. 
Then
$$
\e_{\bari} = -\e_i
$$ 
for all $1 \leq i \leq  n$
and the characters
$\e_1,  \dots, \e_n$ 
form a basis of the character group of $T_{2n}$.
\end{numar}

\begin{numar} \label{nr:liespd} 
The Lie algebra of the symplectic group identifies with the subalgebra
of $\liegld$
$$
\liespd =\{ X\in \liegld \mid  \transp{X} J + J X = 0\}.
$$ 
We have
$X=(x_{ij}) \in \liespd$ 
if and only if
$$
x_{ij}=-x_{\bar{\jmath}\bar{\imath}},\;
x_{i\bar{\jmath}}=x_{j\bar{\imath}},\;
x_{\bar{\imath}j}=x_{\bar{\jmath}i}
$$
for all $1\leq i,j \leq n.$
Let $E_{ij}$ by the elementary matrix with a 1 in the $i$ row and $j$
column, and denote
\begin{align*} 
X_{ij} &=E_{ij}- E_{\barj \bari} \text{ for }1\leq i,j \leq n \text{ and}\\
X_{i\barj} &=E_{i\barj}+E_{j\bari},\;\;
X_{\bari j}=E_{\bari j}+E_{\barj i},
\text{ for } 1\leq i\leq j \leq n.
\end{align*}
These elements make a basis of $\liespd$. 

The elements $X_{ii}$, $1\leq i \leq n$, form a basis of the
Cartan subalgebra $\lietd$ of diagonal matrices in $\liespd$.
The elements $X_{ij}, X_{i\barj}$, $1\leq i \leq j \leq n$, form a
basis of the Borel subalgebra $\liebd$ of $\liespd$ of upper triangular matrices.
The roots of $\liespd$ are
\begin{align*} 
\pm &(\e_i-\e_j), \quad 1\leq i < j \leq n \\
\pm &(\e_i+\e_j), \quad 1\leq i \leq j \leq n
\end{align*}
and for 
$1\leq k \neq \ell \leq 2n$, 
the root space of $\liespd$ corresponding to the root
$\e_k-\e_{\ell}$
is
$\C X_{k\ell}$.  
In particular, the roots of the Borel subalgebra $\liebd$, which we
choose as positive roots, are 
$\e_i-\e_j$, $1\leq i < j \leq n$, and  $\e_i+\e_j$, $1\leq i\leq j \leq n$.
\end{numar}

\begin{numar}     \label{nr:wspd} 
The Weyl group 
$\wspd = N_{\spd}(T_{2n})/T_{2n} $  
of the symplectic group 
$\spd$ 
is isomorphic to the group of 
linear transformations of $F$ which permute the vectors $e_i$ of the
symplectic basis and commute with the matrix $J$.
It identifies with
the group of permutations of the set
$\{1,\dots, 2n\}$ which commute 
with the involution $i\mapsto \bar{\imath}$
\begin{equation} \label{eq:def_wspd}  
\wspd = \{w \in S_{2n} \mid w( \bar{\imath} ) = \b{w(i)} \},
\end{equation}
the correspondence being given by
$$w.e_i=e_{w(i)},\qquad 1\leq i \leq 2n, \; w\in\wspd.$$ 
 Then the action of $\wspd$ on the character group of the maximal torus
 is given by the rule
 $$w.\e_i=\e_{w(i)},\qquad 1\leq i \leq n, \; w\in\wspd.$$ 
Any element  
$w \in \wspd$ is determined by its values on the set
$\{1,\dots,n\}$.
To represent $w$
we use the one-line notation 
$w(1) w(2) \dots w(n)$ with the convention to write $\bari$ instead of
$j= \bari$ when $j \geq n+1$. For example, for $n=4$,  $1674$ gets
written $1\bar3\bar24$.
Given that 
$\e_{\bar{\imath} }=-\e_i$, the bars in the one-line notation
correspond to sign changes in the action on the characters, and
$\wspd$ identifies in this way to the \emph{signed permutations} of
the  $\e_1, \dots, \e_n$.

The group $\wspd$ is a  Coxeter group generated by the reflections
corresponding to the simple roots. The \emph{length} $\ell(w)$ 
of an element $w$ is the minimal number of simple reflections required to
express $w$ as a product of reflections.
The longest element of  $\wspd$ is $\bar1\bar2\dots\bar{n}$ and its
length is $n^2$. 

The \emph{Bruhat order} on the Weyl group $\wspd$ is defined by  
$w\leq w'$ if and only if there is a chain
$$w=w_1\to w_2\to \dots \to w_d=w'$$
such that, for all 
$1\leq i\leq d-1$,  
$\ell(w_{i+1})=\ell(w_i)+1$ and
$w_{i+1}=\sigma _{\alpha } w_i $ 
for a root $\alpha $.
\end{numar}

\begin{numar} \textbf{Symplectic flag manifolds.}
\label{nr:gokd_fod} 
For $1\leq d_1 < d_2 < \dots < d_r \leq n$ a sequence of integers, 
 denote by $\fo(d_1, \dots, d_r, F)$ the projective variety
 of flags of isotropic subspaces
 $$
   \{(V_{d_1} \subset\dots\subset V_{d_r} \subset F ) \mid \dim V_{d_i} =d_i,\, V_{d_i}  
        \text{ isotropic for all } 1\leq i\leq r\}.
 $$
 The symplectic group acts naturally on this variety and the action is
 transitive, ie this variety is a homogeneous space $\spd/P$ with $P$ a parabolic
 subgroup.
 Conversely, any homogeneous space $\spd/P$ with $P$ a parabolic
 subgroup is isomorphic to one of the varieties  $\fo(d_1, \dots, d_r, F)$.
 We call them \emph{symplectic flag manifolds}.

When $r=1$ we will write $\go(k,F)$ instead of $\fo(k,F)$ and call the
varieties $\go(k,F)$
\emph{symplectic grassmannians}. These are subvarieties of the usual
grassmannians $G(k,F)$ and are isomorphic to 
 the homogeneous spaces 
$\spd/P$ with $P$ a maximal parabolic subgroup.
When $r=n$ we denote
$\fo(F)=\fo(1,\dots,n, F)$
and, if there is no risk of confusion, simply call $\fo(F)$ the
\emph{symplectic flag manifold}. This is isomorphic to the homogeneous
space $\spd/B$ with $B$ a Borel subgroup.
Since
the isomorphism class of these varieties depends only on the dimension of $F$, 
we will also write $\gokd$ and $\fod$ 
instead of $\go(k,F)$ and, respectively, $\fo(F)$.
\end{numar}

\begin{numar} \textbf{Schubert cells in symplectic grassmannians.} 
Let's fix in  $F$ the flag
$E_{\bullet}=(E_1\subset E_2\subset\dots\subset E_{2n})$
generated by the symplectic basis
$ \{e_1,\dots,e_{2n}  \} $ ie such that
$E_i=\Spc{e_1, \dots, e_{i}}$ for all $1\leq i \leq 2n$ (we also say that
$ \{e_1,\dots,e_{2n}  \} $ is an \emph{adapted basis} for the flag 
$E_{\bullet}$). 
\end{numar}

\begin{numar} \label{nr:schubert_spd} 
The incidence of a $k$-dimensional subspace $V\subset F$ 
with the flag
 $E_{\bullet}$
is given by the integers
$r_i=\dim (V\cap E_i)$. 
Let $i_1<i_2<\dots<i_k$ be the positions of the $k$ jumps in the sequence
$$0=r_0\leq r_1\leq\cdots\leq r_{2n-1}\leq r_{2n}=k,$$
that is 
$$i_{\alpha }=\min \{i\mid 1\leq i\leq 2n,\;r_i=\alpha
     \},\qquad\text{for }  1\leq \alpha \leq k.$$
We say that the multi-index 
$I=(1\leq i_1<i_2<\dots<i_k\leq2n)$
is the \emph{incidence type}
of $V$ with the flag $E_{\blt}$. 

Not all indices appear as incidence types of isotropic subspaces. We call
those who do \emph{admissible}. These are precisely those
$I=(1\leq i_1<\dots<i_k\leq2n)$  for which 
$i_{\alpha} \neq \b{i_{\beta}}$ for all 
$1 \leq \alpha , \beta \leq k$, that is, for each $1\leq i \leq n$, at most
one of $i$ or $\bari$ appears in $I$.

\begin{notation} \label{notation:iod} 
We denote $\indice{k}{2n}$ 
the set of multi-indices
$(1\leq i_1<i_2<\dots<i_k\leq 2n)$ 
and 
$\iod$ the subset of $\indice{k}{2n}$ of admissible indices.
\end{notation}
                                
The \emph{Schubert cell} in the symplectic grassmannian $\gokd$
associated (with respect to the flag $E_{\bullet}$) to an admissible multi-index
$I\in \iod$ 
is the set of isotropic subspaces of dimension $k$ whose incidence
type is $I$, ie
\begin{equation*} 
 C_I 
   =\bigl\{V\in \go(k,2n) \,\big|\,
   \dim(V\cap E_i)=\alpha,\;\; i_{\alpha}\leq i<i_{\alpha+1},\, 1\leq\alpha\leq k\bigr\}, 
\end{equation*}
where we put $i_{k+1}=2n+1$. 
If $B_{2n}$ denotes the Borel subgroup of $\spd$ which fixes the flag
$E_{\bullet}$, then the Schubert cells, which are clearly
$B_{2n}$-stable, are actually $B_{2n}$-orbits. 
The symplectic grassmannian $\go(k,F)$ decomposes as a disjoint union
$$G_{\omega }(k,F)=\bigcup_{I\in\iod} C_I.$$ 
\end{numar}

\begin{numar} \label{nr:echelon_gokd} 
The Schubert cell $C_I$ is the set of isotropic subspaces of dimension
$k$ which, 
with respect to the symplectic basis $ \{e_1,\dots, e_{2n}\}$, 
can be generated by the rows of a row-echelon matrix of type $I$, that
is a $k\times 2n$ matrix $(a_{\alpha i})$ for which each entry  
$a_{\alpha i_{\alpha } }$ is non zero and all the entries to 
its right
are zero. 
Each $V\in \C_I$ is generated by the rows of a \emph{unique} such
matrix if we require additionally that the entries 
$a_{\alpha   i_{\alpha } }$ be 1 and the entries below an entry
$a_{\alpha i_{\alpha } }$ be zero.
The remaining entries are free entries except for those below an
entry \emph{opposed} to a 1 
(ie below an entry $a_{\alpha \b{i_{\alpha }} }$) which are
determined. The free entries define an isomorphism between $C_I$ and
the affine space $\bbA^{\dim C_I}$. 

For example, in the symplectic grassmannian $\go(3,8)$ the Schubert
cell defined by the multi-index $(4,6,8)=(4,\bar3,\bar1)$ is given by
the row-space of the matrix
$$
\begin{pmatrix} 
* & * & *     & 1   & 0    & 0 & 0 & 0 \\ 
* & * & *     & 0   &\blt & 1 & 0 & 0 \\ 
* & * &\blt  &  0  &\blt & 0 & * & 1
\end{pmatrix}
$$
We have marked by $*$ the free entries and by $\blt$ those which are
determined.
\end{numar}



\begin{numar} 
On the set $\idn$ we consider the order relation for which $I \leq J$
if and only if
$i_\alpha  \leq j_\alpha$ for all $1 \leq \alpha \leq k$.
We consider the induced order on the subset 
$\iod\subset \idn$ of admissible indices.

Let $P_k$ be the parabolic subgroup of $\spd$ which stabilizes the
subspace $E_k=\Spc{e_1,\dots, e_k}$, so that $\go(k,F)=\spd/P_k$. Then
$P_k$ is the maximal parabolic subgroup which misses the simple root
$\alpha _k=\e_k-\e_{k+1}$. The Weyl group $W(P_k)$ of $P_k$ is
isomorphic to $S_k \times W( \Sp_{2(n-k)} )$.   
The coset space $W(\spd)/W(P_k)$ is in bijection with the $T_{2n}$-fixed
points in $\spd/P_k$, and therefore with $\iod$. Explicitly, 
this bijection associates to 
the class of $w\in\wspd$ the
multi-index 
$\{ w(1),\dots,w(k) \} \!\! \uparrow $ obtained by arranging in increasing order the 
elements of the set $\{ w(1),\dots,w(k) \}$.
Via this bijection, the order on $\iod$ corresponds 
with the order on $W(\spd)/W(P_k)$ inherited from the Bruhat order on $\wspd$, 
and so with the Bruhat order on the set of Schubert cells. That is, for all
$I,J\in \iod$, 
$$C_I \subset \b{C_J} \iff I\leq J.$$
\end{numar}

\begin{numar} \label{nr:schubert_vars_gokd}
The \emph{Schubert subvarieties}
of the symplectic grassmannian $\go(k,F)$ are the closures of the 
Schubert cells.
Since the Bruhat order of the Schubert cells agrees with the order on $\iod$, 
the Schubert subvariety $X_I=\b{C_I}$ is given by
$$X_I=\bigcup_{J\leq I} C_J.$$ 
Therefore, in terms of the incidence with the flag $E_{\bullet}$, the Schubert
subvariety $X_I$ is characterized by 
$$X_I=\{V\in\go(k,F) \mid \dim(V \cap E_{i_\alpha } ) \geq \alpha, \;\; 
           1 \leq \alpha  \leq k\}.$$ 
\end{numar}

\begin{numar} \textbf{Schubert cells in the symplectic flag manifold.}
\label{nr:schubert_fod} 
For a signed permutation $w\in\wspd$ the \emph{rank function}
$r_w$ is the function which to a pair 
$(i,j)\in\{1,\dots,n\}\times \{1,\dots, 2n\}$ associates the number
$$r_w(i,j)=\card\{k \mid  k\leq i, w(k)\leq j\}.$$
The integer $r_w(i,j)$ is therefore  
the number of points in the graph 
$\{(k,w(k)) \mid 1\leq k\leq n\}$ of $w$,
situated in the rectangle $\{1,\dots,i\}\times \{1,\dots, j\}$.

The \emph{Schubert cell} $C_w$ of the symplectic flag manifold $\fo(F)$
associated (with respect to the flag $E_{\bullet}$) to a signed permutation 
$w\in\wspd$ 
is the set of flags of isotropic subspaces 
$V_{\blt}=(V_1\subset\dots\subset V_n)$
which  verify
$$\dim(V_i \cap E_j) = r_w(i,j)$$
for all $1\leq i \leq n$ and 
all $1 \leq j \leq 2n$.
The Schubert cells coincide with the orbits of the Borel subgroup $B_{2n}$ of
$\spd$ which fixes the flag $E_{\bullet}$.
Note that we have $V_{\blt}\in C_w $ if and only if
$V_i \in C_{\{w(1),\dots,w(i)\}\uparrow}$ for all $1\leq i \leq n$.
\end{numar}

\begin{numar}   \label{nr:schubert_matrice_fod} 
The Schubert cell
$C_w$ 
is the set of flags of isotropic subspaces which, with respect to the
symplectic basis $\{e_1,\dots,e_{2n}\}$, admit as adapted basis the rows of
a $n \times 2n$ matrix  $(a_{i j})$ for which
$a_{iw(i)}=1$ and $a_{ij}=0$ for all 
$1\leq i \leq n$ and $1 \leq j \leq 2n$,  $j>w(i)$.
Note that the entries which are preassigned the value 1 correspond to the points
of the graph of $w$.
\end{numar}



\begin{numar} \label{nr:schubert_vars_fodd}
The Bruhat order on  $\wspd$ agrees with the Bruhat order on the set of Schubert cells, that
is, for all $w, w' \in \wspd$, 
$$C_w \subset \b{C_{w'} } \iff w \leq w'.$$
The \emph{Schubert subvarieties} of $\fo(F)$ are the closures of the
Schubert cells. The Schubert subvariety $X_w=\b{C_w}$ is then given by
$$X_w=\bigcup_{w'\leq w} C_{w'}.$$ 
The Bruhat order on $\wspd$ and the order on the sets of admissible indices $\iod$ are related
by  \emph{Proctor's criterion} (cf. \cite{Proctor3}) which states that,
for $w, w' \in \wspd$,
\label{eq:crit_proctor} 
$w\leq w'$ if and only if 
$$\{w(1),\dots,w(i)\} \!\! \uparrow \; \leq  \{w'(1),\dots,w'(i)\} \!\!\uparrow$$
for all  $1\leq i \leq n$.
Geometrically this means that, for any $V_{\bullet}=(V_1\subset\dots\subset V_n)$
in $\fo(F)$, 
$V_{\bullet} \in X_w$ if and only if  
$V_i \in X_{\{w(1),\dots,w(i)\}\uparrow}$ for all $1\leq i \leq n$.
In terms of the incidence with the flag $E_{\blt}$, the Schubert subvariety $X_w$ is the
set of flags of isotropic subspaces  
$V_{\blt}=(V_1\subset\dots\subset V_n)$
which  verify
$$\dim(V_i \cap E_j) \geq r_w(i,j)$$
for all $1\leq i \leq n$ and 
all $1 \leq j \leq 2n$.
\end{numar}

\section{The odd symplectic group}


 Let $E$ be a complex vector space of dimension $2n+1$ and $\omega \in
 \Lambda^2E^*$ a generic skew-form on $E$. 
More precisely, we assume
that $\omega  $
is of maximal rank, ie $\rank \omega  =2n$, since the skew-forms of maximal rank  
form an open orbit under the action of the general linear group
$\GL(E)$ on $ \Lambda^2E^*$.

\begin{defin} 
The \emph{odd symplectic group}
$\spdd$ is the stabilizer in
$\GL(E)$ of the skew-form $\omega $.
\end{defin}

Being an isotropy group for an algebraic action, $\spdd$ is a closed
subgroup of $\GL(E)$.
Let $R$ be the one dimensional kernel of $\omega $ with a generator $e_0$, 
choose a supplement
 $F$, so that the restriction of $\omega $ to $F$ is non degenerate,
 and take $ \{e_1,\dots,e_{2n}\} $ to be a symplectic basis of $F$.
 Then in the basis $ \{e_0, \dots, e_{2n}  \} $ the skew-form $\omega$
is given by
$$
\omega (e_i,e_j) = \delta _{i, 2n+1-j}, 
       \quad \text{for all }0\leq i,j \leq 2n.
$$
We will call such a basis an \emph{odd symplectic basis}.
We continue to use the notation $\bari=2n+1-i$ 
 except that now we
 consider $\bar0$ to be not defined.

With respect to the decomposition $R \oplus F = E$ the odd symplectic
group is the group of the matrices of the form
\begin{equation}  \label{eq:spdd_matricial}  
\begin{pmatrix} \lambda & \ell \\ 0 & a 
\end{pmatrix}  
\end{equation}
where $\lambda \in\cst$, $a\in\spd$, $\ell\in\C^{2n}$. 
From now on, we fix the odd symplectic basis $ \{e_0,e_1, \dots,
e_{2n}   \}$
and identify $\spdd$ with 
this group of matrices.

The symplectic
group $\spd$ embeds in $\spdd$ as those matrices \eqref{eq:spdd_matricial} 
with $\lambda =1$ and $\ell=0$.
Denote $U$ the subgroup of $\spdd$ given by the matrices 
 \eqref{eq:spdd_matricial} with $\lambda =1$ and $a=1$.
It is isomorphic to the additive group $(\cd,+)$ and we easily check
that it is a normal subgroup in $\spdd$. 
Therefore the odd symplectic group $\spdd$ is equal to the 
semi-direct product 
$$
(\C^*\times\Sp_{2n})\ltimes U
$$ 
of his two subgroups $\cst\times\spd$ and $U$.
It follows that $\spdd$ is connected and of dimension $(n+1)(2n+1)$.
With the action of $\cst \times \spd$ by interior automorphisms, $U$ is isomorphic
to the dual $F^*$ of the standard representation of $\spd$ (where
$\cst$ acts by homotheties).

Since $U$ is a normal unipotent subgroup, it follows that $\spdd$ is
not reductive.
Actually $U$ is the unipotent radical $R_u$ of $\spdd$, since on the
one hand we have $U\subset R_u$ and on the other hand the quotient
$\spdd/U\simeq \csp$
is reductive so the image of $R_u$ in this quotient is trivial.

We describe now the Borel subgroups and the maximal tori of $\spdd$.

\begin{prop} 
Let $B$ be the subgroup of $\spdd$ 
of upper triangular matrices in the
odd symplectic basis $ \{e_0,e_1, \dots, e_{2n}   \}$. Then $B$ is a Borel subgroup.
\end{prop}

\begin{proof} 
$B$ is solvable since it is a subgroup of the subgroup of upper
triangular matrices in $\gldd$. It is given by the matrices of the form
$$
\begin{pmatrix} \lambda & \ell\\ 0 & a 
\end{pmatrix}   
$$
with $a\in\spd$ upper triangular, so 
$B$ is the semi-direct product
$B=(\cst \times B_{2n})  \ltimes U$
where $B_{2n} \subset \spd $
is the Borel subgroup of  $\spd$ preserving the flag
generated by the symplectic basis $ \{e_1, \dots, e_{2n}   \}$.
In particular $B$ is connected and therefore it is contained in a Borel
subgroup $\tilde{B}$.
Since $\tilde{B}$ contains $U$, 
it is the semi-direct product
$\tilde{B} =[\tilde{B} \cap(\csp)]\ltimes U$.
But $\tilde{B} \cap(\csp)$ is connected 
(since $\tilde{B}$ and $U$ are)
and solvable and contains 
$\C^*\times B_{2n}$. It follows that 
$\tilde{B} \cap(\csp)=\C^*\times B_{2n}$, 
and this means that $\tilde{B} =B$.
\end{proof}

\begin{prop}    \label{prop:tore_max_spdd} 
Let $T$ be the subgroup of $\spdd$ of diagonal matrices in the odd
symplectic basis $ \{e_0,e_1, \dots, e_{2n}   \}$. Then $T$ is a maximal torus.
\end{prop}

\begin{proof}
$T$ is a torus since it is a group of diagonal matrices.
We have $T=\C^*\times T_{2n} $ where $T_{2n}  \subset \spd$
is the maximal torus of $\spd$ of diagonal matrices in the symplectic 
basis  $ \{e_1, \dots, e_{2n}   \}$. It suffices to show that $T$ is
equal to its centralizer.
An element
$$
\begin{pmatrix} \lambda & \ell\\ 0 & a 
\end{pmatrix}
$$
of $\spdd$ centralizes $T$ if and only if
$ad=da$ and $\ell d=t\ell$ for all $(t,d)\in \C^*\times T_{2n} $, 
ie $a\in C_{\spd}(T_{2n} ) $ and $\ell=0$.
Since $C_{\spd}(T_{2n} )= T_{2n} $, we are done.
\end{proof}

Note that the elements of the maximal torus $T$ are the diagonal
matrices
$$\diag(t_0,t _1,\dots,t_n,t_n^{-1}, \dots, t_1^{-1}), \quad t_0, \dots, t_n\in\C^*.$$

\begin{remarca}   \label{rem:gr_weyl_spdd} 
The last proof shows that 
$C_{\spdd}(T)=T$ so the Weyl group of $\spdd$ is the quotient
$N_{\spdd}(T)/T$. 
This is actually
isomorphic to
the Weyl group of the symplectic group $\spd$.
This follows from the fact that 
$N_{\spdd}(T)=\cst\times N_{\spd}(T_{2n} )$,
which we can easily check.
\end{remarca}

For later use, we record the following :

\begin{prop} \label{prop:centre_spdd}                 
The center of $\spdd$ is $ \{\pm \Id \}$.
\end{prop}


\begin{numar} \textbf{Relating $\spdd$ and $\spddd$.} \label{nr:tldo} 
Embed $E$ as a hyperplane in a complex vector space 
$\tE$ of dimension $2n+2$.
Then the odd symplectic form 
$\omega $ extends to a symplectic form $\tldo$ on $\tE$. Indeed it
suffices to take a vector $e_{\bar0}$ in $\tE$ not belonging to $E$
and to define $\tldo$ in such a way that the basis
$\{e_0,e_1, \dots, e_n, e_{\bar{n}}, \dots, e_{\bar1},  e_{\bar0}  \}$
be symplectic. 
Denote simply by $\spddd$ the symplectic group which
fixes $\tldo$.
Let $P$ be the parabolic subgroup of $\spddd$ which preserves the line
$R$. Then $P$ also preserves $E$ which is the orthogonal of $R$, and
for any element $g \in P$ the restriction $\res{g}{E}$ preserves
$\omega $, ie is in $\spdd$.  
\end{numar}

\begin{prop} \label{prop:P->Sp} 
The morphism $P \to \spdd$ 
given by
$g\mapsto \res{g}{E}$ is  surjective.
\end{prop}

\begin{proof}
In the symplectic basis 
$\{e_0,e_1, \dots, e_n, e_{\bar{n}}, \dots, e_{\bar1},  e_{\bar0}  \}$
an element of  $P$ is a matrix of the form 
\begin{equation} \label{eq:matriceP}
\begin{pmatrix} 
\lambda  &  \ell  &  \nu \\
0        &  a  &  c   \\
0        &  0  &  \mu 
\end{pmatrix}
\end{equation}
with  $a\in\spd$, $\lambda, \mu  \in\C^*$, $\nu\in\C$ and
$\ell,c\in\C^{2n}$ a row and, respectively, column
vector.
The condition that this be an element of $\spdd$ is that the columns
make a symplectic basis. 
Since $\res{g}{E}\in\spdd$, 
we only need to look at the conditions involving the last column and
this gives
$\lambda \mu =1$ and, in matrix form,
\begin{equation} \label{eq:c}
\mu \ell+(c_{2n}, \dots, c_{n+1}, -c_n, \dots, -c_1)a=0.
\end{equation}
This shows that given  any
$(\lambda,\ell,a)$ there is a $c$ such that 
\eqref{eq:matriceP} belongs  to $\spddd$ and we are done.
\end{proof}

\begin{remarca}
In the proof we see that giving $(\lambda,\ell,a)$  uniquely
determines $c$ while $\nu $ is arbitrary. We can show that there is no
way of choosing a $\nu $ for each $(\lambda,\ell,a)$ such that the
resulting injection $\spdd \to P$ be a morphism of algebraic groups,
ie the morphism $P \to \spdd$ has no section.
Actually we can show that $\spdd$ cannot be embedded as an algebraic
subgroup in $\spddd$. For a proof of these statements, see \cite[ch. 1]{teza}. 
\end{remarca}

\begin{remarca}
Through the morphism $P\to \spdd$, the Borel subgroup of $\spddd$ of
upper triangular matrices in the symplectic basis
$\{e_0,e_1, \dots, e_n, e_{\bar{n}}, \dots, e_{\bar1},  e_{\bar0}  \}$
surjects onto the Borel subgroup of $\spdd$ of upper triangular
matrices in the odd symplectic basis 
$\{e_0,e_1, \dots, e_n, e_{\bar{n}}, \dots, e_{\bar1}  \}$.
Similarly, the morphism $P \to \spdd$ restricts to an isomorphism
between the maximal torus of $\spddd$ of diagonal matrices
in the basis
$\{e_0,e_1, \dots, e_n, e_{\bar{n}}, \dots, e_{\bar1},  e_{\bar0}  \}$
and the maximal torus of $\spdd$ of diagonal matrices
in the basis 
$\{e_0,e_1, \dots, e_n, e_{\bar{n}}, \dots, e_{\bar1} \}$.
 \end{remarca}

\begin{remarca} 
From \eqref{eq:c} it follows that the kernel of the morphism
$P\to\spdd$ is the unipotent group of dimension one whose elements are
the matrices of the form
 \begin{equation}\label{eq:ker_restr}  
 \begin{pmatrix}
 1 & 0 & \nu \\
 0 & 1 & 0   \\
 0 & 0 & 1
 \end{pmatrix} 
 \end{equation}
with $\nu \in \C$.
\end{remarca}

\begin{numar} \textbf{The intermediate odd symplectic group.} \label{nr:intermediate} 
Since $\spdd$ does not embed in $\spddd$ as an algebraic subgroup, it
is desirable to have an alternative odd symplectic group which sits
between $\spd$ and $\spddd$. 
Such a group has been introduced by
Gelfand and Zelevinski  in \cite{GelfZelev}.
They notice that for any complex vector space $V$, the general linear
group
$\GL(V)$ has an open orbit in the space
$V^* \oplus \Lambda ^2V^*$,       
which is 
$$\{(\ell,\theta ) \mid \ell \neq 0,\;\theta \text{ of maximal rank and} \; \ker \ell \cap
        \ker \theta =(0)\},$$
and they define the group $G(V)$ to be the isotropy group of this orbit.
It is easy to check that if $V$ is of dimension $2n+1$ then $G(V)$ is
isomorphic to $\spd$, while if $V$ is of dimension $2n+2$ then
$G(V)$ is isomorphic to the subgroup of $\spddd$ which fixes a
non-zero element of $V^*$ (or, equivalently, of $V$).
In the latter case, we will use the following designation for $G(V)$:

\begin{defin} \label{def:tspdd} 
The \emph{intermediate odd symplectic group} 
$\tspdd$ is 
the stabilizer in $\spddd$ of a non-zero element of  $\cddd$.
\end{defin}

Since $\tspdd$ is well defined up to a conjugation, we can assume that
the element it fixes is the first vector $e_0$ of the symplectic basis
$\{e_0,e_1, \dots, e_n, e_{\bar{n}}, \dots, e_{\bar1},  e_{\bar0}  \}$
so that $\tspdd$ identifies with the subgroup of the parabolic $P$
considered above given by the matrices \eqref{eq:matriceP} with 
$\lambda = \mu = 1$. Its image via the morphism $P \to \spdd$ is the
subgroup $\sspdd$ of elements of determinant 1 in $\spdd$.
\end{numar}

\section{Odd symplectic flag manifolds}

We introduce here our main objects of study.
For a sequence of integers $1\leq d_1 < \dots < d_r \leq n+1$, denote by
$\fo(d_1,\dots,d_r,E)$ the variety of flags
 $$
   \{(V_{d_1} \subset\dots\subset V_{d_r} \subset E ) \mid \dim V_{d_i} =d_i,\, V_{d_i}  
        \text{ isotropic for all } 1\leq i\leq r\}.
 $$
These are the \emph{odd symplectic flag manifolds}.

Let us first note that the condition that the integers $d_i$ be at
most $n+1$ comes from the fact that $n+1$ is the maximal dimension of an
isotropic subspace of $E$. Actually, a maximal isotropic subspace of
$E$ is always of dimension $n+1$ since it contains the kernel $R$ and
its image in $E/R$ is a maximal isotropic subspace.

When $r=1$ we will write $\go(k,E)$ instead of $\fo(k,E)$ 
and call
these varieties
\emph{odd symplectic grassmannians}. These are simply the
grassmannians of isotropic subspaces of $E$.
When $r=n+1$ we denote $\fo(E)=\fo(1,\dots, n+1,E)$ and call this, if
there is no risk of confusion, the \emph{odd symplectic flag manifold}.
This is the variety of maximal flags of isotropic subspaces of $E$.

Since
the isomorphism class of the variety
$\fo(d_1,\dots,d_r,E)$ depends only on the dimension of $E$, we will 
also denote this variety by $\fo(d_1,\dots,d_r; 2n+1)$ (and correspondingly, 
we will also use the notations $\gokdd$ and $\fodd$).

The variety $\fo(d_1,\dots,d_r,E)$ is a closed subvariety of
$\bbF(d_1,\dots,d_r,E)$, the type $\ssA$ flag manifold
defined by forgetting the isotropy
condition, therefore the odd symplectic flag manifolds are projective
varieties. The odd symplectic group acts naturally on the varieties
$\fo(d_1,\dots,d_r,E)$, the action being defined by restricting to 
$\fo(d_1,\dots,d_r,E)$ the usual action of $\GL(E)$ on the flag
manifold $\bbF(d_1,\dots,d_r,E)$. The major difference between this
situation and the one we have in the symplectic setting is that this
action is not transitive. The reason, as explained in the
introduction, is that the kernel $R$ is fixed by $\spdd$ and so flags
which do not have the same incidence with $R$ cannot be
conjugated by $\spdd$ (we still have to prove that these different
incidence types actually occur, but this will become clear below).
We will show that actually these incidence conditions suffice to describe
the orbits.

Being non homogeneous, it is no longer granted that the odd symplectic
flag varieties
$\fo(d_1,\dots,d_r,E)$ are smooth. 
It turns out that they actually are:


\begin{prop} 
The odd symplectic flag manifold $\fo(d_1,\dots,d_r,E)$ is a smooth subvariety of 
codimension $\frac12 d_r(d_r-1)$ in the flag manifold 
$\bbF(d_1,\dots,d_r,E)$.
\end{prop}

\begin{proof}
Let $T=T_{d_r}$ denote the highest rank tautological bundle on the flag manifold   
$\bbF(d_1,\dots,d_r,E)$. The fiber of $T$ at a point 
$V_{\bullet}= (V_{d_1} \subset\dots\subset V_{d_r})$
is $V_{d_r}$. 
Any skew-form $\theta \in \Lambda ^2E^*$ can be seen as a global section of
$\Lambda ^2T^*$, whose value at a point
$V_{\bullet}= (V_{d_1} \subset\dots\subset V_{d_r})$
is the restriction $\res{\theta}{V_{d_r} }\in\Lambda ^2 V_{d_r} ^*$.
In this way the section defined by 
$\theta \in \Lambda ^2E^*$ 
vanishes at a point $V_{\bullet}$ if and only if $V_{d_r}$ is isotropic with 
respect to $\theta $. Therefore the zero locus of our odd symplectic form   
$\omega $ seen as a section of  $\Lambda ^2T^*$
is $\fo(d_1,\dots,d_r,E)$.
Now, the vector bundle $\Lambda ^2T^*$ is generated
by its global sections which come from $\Lambda ^2E^*$
since any skew-form defined on a subspace $V_{d_r}$ extends to $E$.
The section defined by the odd symplectic form $\omega $ is generic among the 
sections which come from  $\Lambda ^2E^*$, so to conclude by Bertini's
theorem it suffices 
to show that $\omega $ vanishes in at least a point of 
$\bbF(d_1,\dots,d_r,E)$. For this it is enough to take 
the point $(E_{d_1} \subset\dots\subset E_{d_r})$ where $E_{\bullet}$ is the
complete flag generated by the odd symplectic basis 
$\{ e_0, e_1, \dots, e_{2n} \}$.
\end{proof}

\begin{numar}  \textbf{$\spdd$-orbits in the odd symplectic grassmannians.}
\label{nr:orbits_in_gokdd} 
Any isotropic subspace $V\subset E$ of dimension $n+1$ contains the kernel $R$
and is of the form $V=R \oplus W$ with $W\subset F$ isotropic of dimension $n$, 
so the odd symplectic grassmannian $\go(n+1,E)$ is isomorphic to the 
symplectic grassmannian $\go(n,F)$. 
Therefore $\go(n+1,E)$ is already homogeneous under the subgroup $\spd\subset \spdd$.
For the other odd symplectic grassmannians, the $\spdd$-orbits are given by 
the incidence with the kernel $R=\C e_0$, that is we have:

\begin{prop} 
For $1\leq k\leq n$ the odd symplectic group $\spdd$ acts on the
odd symplectic grassmannian $\gokdd$ with two orbits 
\begin{align*} 
X_0 &= \{ V\in\gokdd\mid e_0\in V\} \\
X_1 &= \{V\in\gokdd\mid e_0\notin V\}.
\end{align*}
The closed orbit $X_0$ is isomorphic to the symplectic grassmannian
$\go(k-1,2n)$ and the open orbit $X_1$ is isomorphic to the total space
of the dual of the tautological bundle on the symplectic
grassmannian $\gokd$.
\end{prop}

\begin{proof}
Both $X_0$ and $X_1$ are obviously stable under $\spdd$.
Any $V\in X_0$ is of the form $V=R \oplus W$ with $W=V \cap F$ isotropic of 
dimension $k-1$, so $X_0$ is isomorphic to the symplectic grassmannian
$\go(k-1,F)$ and is already an orbit of the subgroup $\spd \subset \spdd$.

Let now $p:E \to F$ be the projection coming from the decomposition 
$E=R \oplus F$. 
For any $V \in X_1$, $p(V)\subset F$ is an isotropic subspace of dimension $k$.
We get in this way a map
$$p:X_1\vers \go(k,F)$$
onto the symplectic grassmannian $\go(k,F)$ which we regard as a subvariety in $\go(k,E)$.
If $V=\Spc{v_1, \dots, v_{k}} \in X_1 $, with $v_i=\alpha _ie_0+v'_i$, 
$v'_i\in F$, then $p(V)=\Spc{v'_1,\dots,v'_k}$.
Therefore $V$ and $p(V)$ will be conjugated by any element
of the unipotent radical $U$ of $\spdd$ which sends
$v'_i\mapsto v'_i+\alpha _i e_0$ for all $1\leq i \leq k$. 
Such an element exist since the $v'_i$ are independent.
Conversely, for any 
$g\in U$ we have $p(g.V)=p(V)$, 
and so the orbits of $U$ in $X_1$ coincide with the fibers
of $p:X_1 \to \go(k,F)$.
In particular, any orbit of $U$ in $X_1$ meets $\go(k,F)\subset X_1$
which is an orbit of $\spd \subset \spdd$, and therefore 
$\spd \ltimes U \subset \spdd$ acts transitively in $X_1$.

Let now $T$ be the tautological bundle on the symplectic grassmannian
$\go(k,F)$. The fiber of the dual $T^*$ at a point $W \in \go(k,F)$
is $W^*$. We define a map 
$T^* \to X_1$ by sending 
an element $\ell\in W^*$ to its graph
$\Gamma _{\ell }\subset W\oplus \C=W\oplus R \subset F\oplus R=E$ 
which is an isotropic subspace of $E$ which does not contain $e_0$ 
and for which 
$p(\Gamma _{\ell})=W$. 
In coordinates, if 
$W=\Spc{w_1,\dots,w_k} $ then
$\Gamma _{\ell}  =\Spc{w_1+\ell(w_1)e_0,\dots,w_k+\ell(w_k)e_0} $.
Conversely, if $V\in X_1$ then $V\subset p(V) \oplus \C e_0$ 
is the graph of a linear map
$p(V)\to \C e_0=\C$, ie an 
element of the fiber of $T^*$ at the point $p(V)$.
\end{proof}

\end{numar}

\begin{numar}  \textbf{$\spdd$-orbits in the odd symplectic flag manifold $\fodd$.}
As for the odd symplectic grassmannians, the $\spdd$-orbits of the odd symplectic 
flag manifold $\fo(E)$ are described by the incidence with the kernel $R=\C e_0$.
More precisely, we have:

\begin{prop} \label{prop:orbites_fodd}
The odd symplectic group  $\spdd$ acts on the odd symplectic flag manifold $\fodd$ 
with $n+1$ orbits
$$
\cO_i=\big\{(V_1\subset \dots\subset V_{n+1})\in\fodd \mid e_0\in V_i, e_0\notin V_{i-1}\big\},
\quad 1\leq i\leq n+1.
$$
Moreover, for all $1\leq i\leq n+1$, the orbit $\cO_i$ is isomorphic to
the total space of the dual $T^*_{i-1}$ of the rank $i-1$ tautological
bundle on the symplectic flag manifold $\fod$.
\end{prop}

\begin{proof} 
The $\cO_i$ are clearly $\spdd$-stable.
Let again $p:E\to F$ be the projection coming from the decomposition $E=R \oplus F$.
For any $V_{\blt}=( V_1\subset \dots\subset V_{n+1})$ in $\cO_i$,
the projection $p$ restricts to isomorphisms $V_j\to p(V_j)$ 
for all $1\leq j\leq i-1$, and we have
$V_j=p(V_j)\oplus R$ with $p(V_j)=V_j\cap F$ 
for $i\leq j\leq n+1$. We also have 
$V_i=V_{i-1}\oplus R$ and $p(V_{i-1})=p(V_i)$. 
Let's denote by $p_i(V_{\blt})$ the flag 
$$
p_i(V_{\blt})=(p(V_1)\subset\dots\subset p(V_{i-1})\subset V_i \subset\dots\subset V_{n+1}).         
$$
We get in this way a map $p_i$ from $\cO_i$ onto the closed subvariety $Y_i \subset \cO_i$
$$
Y_i=\big\{V_{\blt}=(V_1\subset \dots\subset V_{n+1})\in\fo(E) \;\big|\; 
   V_1,\dots,V_{i-1}\subset F, \; e_0\in V_i\big\}.
$$   
The subvariety $Y_i$ is isomorphic to the symplectic flag manifold $\fo(F)$ via the map
which sends a flag $(W_1\subset\dots\subset W_n) \in \fo(F)$ to the flag
$$
(W_1\subset\dots\subset W_{i-1}\subset 
    W_{i-1} \oplus R \subset\dots\subset W_n\oplus R) \in Y_i.
$$
In coordinates, if $\{ v_1', \dots, v_n' \}$ is an adapted basis for the flag in
$\fo(F)$ which corresponds via the isomorphism above to the flag $p_i(V_{\blt})$, 
then 
$$
\{v_1',\dots,v_{i-1}', e_0, v_i',\dots, v_n'\}
$$
is an adapted basis for $p_i(V_{\blt})$, and there are
$a_1,\dots,a_{i-1}\in\C$ such that 
$$
\{v_1'+a_1e_0,\dots,v_{i-1}'+a_{i-1}e_0,e_0,v_i',\dots, v_n'\}
$$
is an adapted basis for the flag 
$V_{\blt}$.
It follows that $V_{\blt}$ and $p_i(V_{\blt})$
will be conjugated by any element of the unipotent radical
$U$ of $\spdd$ which sends $v_j'\mapsto v_j'+a_je_0$
for all $1\leq j \leq i-1$.
Such an element exists since 
$v_1',\dots,v_{i-1}'$ 
are independent.
Conversely, for any
$g\in U$ we have
$p_i(V_{\blt}) =p_i(g\cdot V_{\blt})$, so the orbits of $U$ in $\cO_i$ 
coincide with the fibers of $p_i$.
Therefore any orbit of $U$ in $\cO_i$ meets $Y_i \subset \cO_i$ which
is an orbit of $\spd \subset \spdd$, and so it follows that
$\spd \ltimes U \subset \spdd$ acts transitively on
$\cO_i$.

Let now $T_{i-1}$ be the tautological bundle of rank $i-1$ on the 
symplectic flag manifold $\fo(F)$. The fiber of the dual  
$T^*_{i-1}$ at a point 
$W_{\blt} = (W_1\subset\dots\subset W_n) \in \fo(F)$
is $W_{i-1}^*$.
We define a map $T_{i-1}^* \to \cO_i$ by sending an element  
$\ell\in W_{i-1}^*$ to the flag
$V_{\blt}=(V_1\subset \dots\subset V_{n+1}) \in \cO_i$ 
where: 
$V_{i-1}$ is the graph of $\ell$, 
$\Gamma _{\ell}\subset W_{i-1}\oplus\C=W_{i-1}\oplus R \subset F\oplus R = E$;
for $1 \leq j \leq i-2$, $V_j$ is the graph
$\Gamma _{\res{\ell}{W_j}}\subset \Gamma _{\ell}$ of 
the restriction
$\res{\ell}{W_j}$; 
and 
$V_j=W_{j-1}\oplus R$
for $i\leq j\leq n+1$.
In coordinates, if
$ \{v_1',\dots,v_n' \} $ is an adapted basis for the flag
$W_{\blt}$ then 
$$\{v_1'+\ell(v_1') e_0, \dots, v_{i-1}'+\ell(v_{i-1}') e_0, e_0, v_i', \dots, v_n'\}$$
is an adapted basis for $V_{\blt}$. 
Conversely, for 
$V_{\blt}=(V_1\subset \dots\subset V_{n+1}) \in \cO_i$ 
let $W_{\blt} = (W_1\subset\dots\subset W_n) $
be the flag in $\fo(F)$ corresponding via the isomorphism $\fo(F)\simeq Y_i$
to the flag $p_i(V_{\blt})$, ie 
$$
W_{\blt}=(p(V_1)\subset \dots \subset p(V_{i-1}) \subset p(V_{i+1}) \subset \dots \subset p(V_{n+1})).
$$
Then $V_{i-1}$ is the graph of a linear map 
$\ell:W_{i-1}\to\C e_0=\C$, ie an element of the fiber of $T_{i-1}^*$
at the point $W_{\blt}$, and
we define the inverse map $\cO_i \to T_{i-1}^*$
by sending $V_{\blt}$ to $\ell$.
\end{proof}


\begin{prop} \label{prop:orbit_closure} 
The closures of the $\spdd$-orbits of $\fodd$ are the subvarieties
$$
\overline{\cO_i}=
\big\{(V_1\subset \dots\subset V_{n+1})\in\fodd \mid e_0\in V_i \big\},
\quad 1\leq i\leq n+1
$$
and are smooth.
\end{prop}

\begin{proof} 
The first assertion is clear. Now
let $\pi _i : \foe \to \go(i,E)$ be the  map 
$(V_1\subset\dots\subset V_{n+1}) \mapsto V_i$. The image 
$\pi _i (\b{\cO_i})$ is the closed $\spdd$-orbit $X_0$ in the odd symplectic
grassmannian $\go(i,E)$. 
For $V_i \in X_0$, the fiber 
$\pi _i^{-1} ( V_i )$ is isomorphic to
$\bbF(V_i) \times \bbF_{\bar\omega} ( V_i^{\bot} / V_i )$   where
$\bar\omega $ is the form induced by $\omega $ on the quotient
$V_i^{\bot} / V_i $. For any $V_i \in X_0$  we have 
$\dim(V_i^{\bot})= 2n+2-i$ so $\pi _i: \cO_i \to X_0$ is a fibration
with fiber $\bbF(i) \times \fo(2(n+1-i))$ over the smooth base
$X_0 \simeq \go(i-1, 2n)$.
\end{proof}


Note that from the second assertion of the proposition
\ref{prop:orbites_fodd} it follows that
$\cO_i$ is of codimension
$n+1-i$ in $\fodd$.
In particular, $\overline{\cO_n}=\fodd\setminus\cO_{n+1}$ is an 
irreducible $\spdd$-stable divisor.


\begin{remarca}
The $\spdd$-orbits in the other (``partial'') odd symplectic flag
manifolds satisfy the obvious analogs of propositions 
\ref{prop:orbites_fodd} and
\ref{prop:orbit_closure}.
\end{remarca}
\end{numar}

\begin{numar} \textbf{Schubert cells in odd symplectic grassmannians.}
\label{nr:cells_in_gokdd} 
Since the odd symplectic grassmannian $\go(k,E)$ is a subvariety of the
usual grassmannian $G(k,E)$, 
it seems natural to try to define 
``Schubert cells'' in $\go(k,E)$ 
by using incidence conditions with respect to some fixed flag, ie to take
Schubert cells in $G(k,E)$ and intersect them with $\go(k,E)$.
And indeed, if the fixed flag is generated by an odd symplectic basis, we can
readily show that the ``Schubert cells'' we obtain in this way can be described
in terms of row echelon matrices, as for the symplectic grassmannians, and so prove, 
for example, that they are isomorphic to affine spaces.
It is then a bit trickier, though, to derive the incidence relations between the cells, 
ie the Bruhat order on the set of cells.
Things 
simplify noticeably once we observe that the odd symplectic grassmannian
$\go(k,E)$ can be identified with a Schubert subvariety of a symplectic
grassmannian $\gokddd$.

Let indeed, as in \ref{nr:tldo},  $\tlde$ be a complex vector space 
of dimension $2n+2$ which contains
$E$ as a hyperplane and $\tldo$ a symplectic form on $\tlde$ which
extends $\omega $. Then, clearly, an element $V \in \go(k,E)$
is nothing but an isotropic subspace of $\tlde$ which is contained in the
hyperplane $E$, that is, $\go(k,E)$ identifies with the subvariety
\begin{equation} \label{eq:schubert_gokdd} 
\{V \in G_{\tldo} (k,\tlde)  \mid V \subset E \}.
\end{equation}
This is a Schubert subvariety in $G_{\tldo} (k,\tlde)$. 
Let $e_{2n+1} \in \tlde$ be such that the basis
$\{e_0, \dots, e_{2n}, e_{2n+1} \}$ is
symplectic and let $E_{\blt}$ be 
the flag generated by this basis.
Then $E=E_{2n+1}$ and
in the notation of \ref{nr:schubert_vars_gokd}, 
the subvariety \eqref{eq:schubert_gokdd} is the Schubert subvariety
$X_{\bar{k}\dots\bar2\bar{1} }$ 
if $k<n+1$, respectively 
$X_{0\bar{n}\dots\bar2\bar{1} }$ 
if $k=n+1$.
We denote by $\mx_{k,n}$ the multi-index of the Schubert variety \eqref{eq:schubert_gokdd}
that is
  \begin{equation*} 
    \mx_{k,n} =
    \begin{cases} 
      (\bar{k},\dots,\bar{1})  & \text{ if } k< n+1 \\ 
      (0,\bar{n},\dots,\bar{1})  & \text{ if } k=n+1.
    \end{cases} 
  \end{equation*}
We will actually suppose that $k$ and $n$ are fixed and slightly abuse notation and
write $\mx$ instead of $\mx_{k,n}$. 

Through its identification with the Schubert variety $X_{\mx}$, 
the odd symplectic grassmannian $\go(k,E)$ 
acquires the cell decomposition of the latter
$$X_{\mx}=\bigcup_{I \leq \mx}C_I.$$
We then define the \emph{Schubert cells} of 
$\go(k,E)$ (with respect to the flag $(E_1 \subset \dots \subset E_{2n+1})$)
to be the Schubert cells of $G_{\tldo}(k,\tlde)$ 
(with respect to the flag
$(E_1 \subset \dots \subset E_{2n+1} \subset E_{2n+2})$)
which are included in $X_{\mx}$.     

We check easily that for $I \in \ioddd$,
$I=(0\leq i_1 < \dots < i_k \leq 2n+1)$, we have $I \leq \mx$ if and only if $i_k \leq 2n$, 
that is if $I$ is a multi-index in the set $\{0,\dots,2n\}$. So let us introduce the following
notation:

\begin{notation}
Let $\idd$ be the set of multi-indices $I=(0\leq i_1 < \dots < i_k \leq 2n)$
and $\iodd$ the subset of \emph{admissible} indices, ie indices $I$ such that, 
for each $1\leq i \leq n$, at most one of $i$ or $\bari$ appears in $I$.
\end{notation}
Then $\iodd$ identifies with the interval $ \{ I \in \ioddd \mid I \leq \mx \}$ so the 
Schubert cells of $\go(k,E)$ are parametrized by $\iodd$. 
We state now the characterization
of the Schubert cells of $\go(k,E)$ by incidence conditions
with respect to the flag 
$(E_1 \subset \dots \subset E_{2n+1})$:

\begin{prop} \label{prop:incidence_cells_gokdd} 
The Schubert cell $C_I$ of $\go(k,E)$ associated 
to an admissible multi-index $I \in \iodd$
is the set of isotropic subspaces of dimension $k$ in $E$ whose incidence type
is $I$, ie
\begin{equation*} 
 C_I 
   =\bigl\{V\in \go(k,2n+1) \,\big|\,
   \dim(V\cap E_i)=\alpha,\;\; i_{\alpha}\leq i<i_{\alpha+1},\, 1\leq\alpha\leq k\bigr\}, 
\end{equation*}
where we put $i_{k+1}=2n+1$. 
\end{prop}

\begin{proof} 
Indeed, these incidence conditions are obtained from those 
defining $C_I$ in $G_{\tldo}(k,\tlde)$ by forgetting the last one (corresponding to $i=2n+1$),
which is superfluous since $i_k \leq 2n$.
\end{proof}
\end{numar}

\begin{numar} \label{no:C_I_matrice}
Similarly to the symplectic case, the Schubert cell $C_I$ in
$\go(k,E)$ is the set of isotropic subspaces 
of dimension $k$ in $E$ which, with respect to the odd symplectic basis 
$\{ e_0, \dots, e_{2n} \}$, 
can be generated by the rows of a row-echelon matrix of type $I$, ie a 
$k\times(2n+1)$ matrix 
$(a_{\alpha i})$ for which each entry
$a_{\alpha i_{\alpha } }  $ is non zero
and all entries to its right are zero.
For each $V \in C_I$ there is a unique such matrix verifying 
further the conditions that the entries
$a_{\alpha i_{\alpha } } $ be 1 and all the entries below be zero.
The remaining entries are free, except the entries below an entry opposed to a 1 
(ie below an entry of the form $a_{\alpha \b{i_{\alpha } } } $) which are
determined. The free entries provide explicitly an isomorphism between $C_I$ 
and the affine space $\bbA^{\dim C_I}$. 

For example, in the odd symplectic grassmannian $\go(3,9)$ the Schubert
cell defined by the multi-index $(4,6,8)=(4,\bar3,\bar1)$ is given by
the row-space of the matrix
$$
\begin{pmatrix} 
* & * & * & *     & 1   & 0    & 0 & 0 & 0 \\ 
* & * & * & *     & 0   &\blt & 1 & 0 & 0 \\ 
* & * & * &\blt  &  0  &\blt & 0 & * & 1
\end{pmatrix}
$$
We have marked by $*$ the free entries and by $\blt$ those which are
determined.
\end{numar}

\begin{numar}
Let $B$ be the Borel subgroup of $\spdd$ of upper triangular matrices
in the odd symplectic basis $\{ e_0, \dots, e_{2n} \}$, ie
the subgroup preserving the flag $(E_1 \subset \dots \subset E_{2n+1})$.
From \ref{prop:incidence_cells_gokdd}, it follows that the Schubert
cells of $\go(k,E)$ are $B$-stable. We actually have:

\begin{prop}
The Schubert cells of $\go(k,E)$ are the orbits of $B$.
\end{prop}

\begin{proof} 
The Schubert cells of $\go(k,E)$ are orbits of the Borel 
subgroup $B_{2n+2}$ of $\spddd$ of upper triangular
matrices in the symplectic basis  
$\{ e_0, \dots, e_{2n+1} \}$.
But $B_{2n+2}$ acts on the Schubert subvariety
$X_{\mx}$ through the restriction morphism $b \mapsto \res{b}{E}$, 
and the image of $B_{2n+2}$ via this morphism is $B$.
\end{proof}
\end{numar}



\begin{numar}
We equip $\iodd$ with the order relation for which $I \leq J$
if and only if $i_{\alpha }  \leq j_{\alpha }$ for all 
$1 \leq i \leq k$. This corresponds with the order on $\ioddd$
when we identify $\iodd$ with the interval 
$ \{ I \in \ioddd \mid I \leq \mx \}$.
It follows that the order on $\iodd$ describes the Bruhat order on the
set of Schubert cells of $\go(k,E)$, that is, for all
$I,J \in \iodd$,
$$C_I \subset \b{C_J} \iff I\leq J.$$
\end{numar} 

\begin{numar}
Define the \emph{Schubert subvarieties} of the odd symplectic
grassmannian $\go(k,E)$ to be the closures of the Schubert cells.
Since the Bruhat order on the set of Schubert cells agrees with the order 
on $\iodd$, the Schubert subvariety $X_I = \b{ C_I }$ verifies 
$$X_I=\bigcup_{J\leq I} C_J.$$ 
Therefore, in terms of the incidence with the flag 
$(E_1 \subset \dots \subset E_{2n+1})$, the Schubert
subvariety $X_I$ is characterized by 
$$
X_I=\{V\in\go(k,E) \mid \dim(V \cap E_{i_\alpha } ) \geq \alpha, \;\; 
        1 \leq \alpha  \leq k \}.
$$
\end{numar}

\begin{numar} 
The Picard group of the odd symplectic grassmannian $\goke$ is
$\Z$. 
Indeed, since $\goke$ has a cellular decomposition, the Picard group
coincides with 
the free abelian group generated by the classes of the 
closures of the codimension 1 cells. In $\goke$ there is only one
codimension 1 cell, ie only one Schubert divisor.
We can easily check that, for $1 \leq k \leq n$,  the Schubert divisor
is 
$$
X_{\b{k+1}\, \b{k-1}\ldots \bar2 \bar1 } 
  = \{  V\in \goke \mid V \cap E_{2n-k+1} \neq 0  \}
$$  
and, for $k=n+1$, it is
$$
X_{0n\b{n-1}\ldots \bar2\bar1 }
  = \{  V\in \go(n+1,E) \mid \dim(V \cap E_{n+1}) \geq 2  \}.
$$
\end{numar}

\begin{numar} \textbf{Schubert cells in the odd symplectic flag manifold $\fodd$.}
\label{sect:cell_sch_fodd} 
In order to define the Schubert cells in the odd symplectic flag manifold $\fodd$, 
we proceed, as in the case of the odd symplectic grassmannians, by first
identifying $\fodd$ to a Schubert subvariety in a symplectic flag manifold $\foddd$.
With the same notations as in \ref{nr:cells_in_gokdd}, a flag
$V_{\blt} \in \fo(E)$ is the same thing as a flag $V_{\blt} \in \fto(\tlde)$
for which $V_{n+1} \subset E$, that is, $\fo(E)$ identifies with the subvariety
$$
\{ (V_1\subset\dots\subset V_{n+1}) \in\bbF_{\tldo}(\tlde) \mid V_{n+1}\subset E\}.
$$
This is a Schubert subvariety in $\fto(\tlde)$. With respect to the flag
$(E_1 \subset \dots \subset E_{2n+2})$, 
and in the notation of \ref{nr:schubert_vars_fodd},
this is the Schubert subvariety 
$X_{\bar{1}\bar{2}\dots\bar{n}0 }$ 
associated to the signed permutation
$ \bar{1}\bar{2}\dots\bar{n}0 \in \wspddd$.

Identifying then $\fo(E)$ with $X_{\wz}$, we define the \emph{Schubert cells}
of $\fo(E)$ (with respect to the flag $(E_1 \subset \dots \subset E_{2n+1})$)
to be the Schubert cells of $\fto(\tlde)$ 
(with respect to the flag $(E_1 \subset \dots \subset E_{2n+1} \subset E_{2n+2})$)
which are included in $X_{\wz}$. 

The Schubert cells of $\fo(E)$ are therefore parametrized by the interval
$$ 
\{ w\in \wspddd\mid w\leq \bar{1}\bar{2}\dots\bar{n}0 \} 
$$
in the Weyl group $\wspddd$ of $\spddd$. We introduce the following notation:
\end{numar}

\begin{notation}
We denote the interval $\{w \in W(\spddd)\mid w\leq \bar{1}\bar{2}\dots\bar{n}0 \} $
by $\wspdd$.
\end{notation}

This notation is somewhat misleading since it might suggest that $\wspdd$ 
is the Weyl group of $\spdd$. Of course, $\wspdd$ is not a group at all, 
and the Weyl group of $\spdd$ coincides with the Weyl group of $\spd$.

\begin{prop} \label{prop:weyl_spdd} 
The elements of $\wspdd$ are precisely those signed permutations 
which in the one-line notation do not contain 
$\bar{0}$, ie those $w$ for which 
$\bar{0}\notin \{w(0),w(1),\dots,w(n)\}$.
\end{prop}

\begin{proof} 
Let $w\in\wspdd$. We have $w \leq \wz$ if and only if the Schubert cell
$C_w$ is contained in $X_{\wz}$ which is equivalent to 
$E_{\blt}^w \in X_{\wz}$, where $E_{\blt}^w \in C_w$ is the 
flag generated by the basis
$\{e_{w(0)}, e_{w(1)}, \dots, e_{w(n)} \}$.
By the definition of the Schubert subvariety $X_{\wz}$ this is equivalent
to $w(i) \in \{0, 1, \dots, 2n\}$ for all $0 \leq i \leq n$ and we are done. 
\end{proof}

\begin{numar} We derive now the characterization of the Schubert cells by incidence
conditions with respect to the flag $(E_1 \subset \dots \subset E_{2n+1})$. 
Let us first note that
if $w \in \wspdd \subset \wspddd$ then, by the last proposition, the rank function 
$r_w$ of $w$ is determined  by its restriction to the rectangle
$\{0,\dots,n\}\times\{0,\dots,2n\}$. 
\end{numar}

\begin{prop} \label{prop:incidence_cond_Cw} 
The Schubert cell $C_w$ of $\fo(E)$ associated to $w\in\wspdd$
is the set of flags of isotropic subspaces 
$V_{\blt}=(V_1\subset\dots\subset V_{n+1})$
which  verify
$$\dim(V_i \cap E_j) = r_w(i,j)$$
for all $0\leq i \leq n$ and 
$0 \leq j \leq 2n$.
\end{prop}

\begin{proof}
These conditions are obtained from those defining $C_w$ in $\fto(\tlde)$ by forgetting
the case $j=2n+1$ for which the conditions are automatically 
verified since $V_{n+1}\subset E$ and
$w \in \wspdd \subset \wspddd$. 
\end{proof}

\begin{numar}
Analogously to the symplectic case, the Schubert cell $C_w$ in $\fo(E)$
is the set of flags of isotropic subspaces which, with respect to the odd symplectic
basis $\{ e_0, \dots, e_{2n} \}$, have an adapted basis given by the rows of
a $(n+1) \times (2n+1)$ matrix $(a_{ij})$ for which
$a_{iw(i)}=1$ and $a_{ij}=0$ for all 
$0\leq i \leq n$ and $0 \leq j \leq 2n$,  $j>w(i)$.
\end{numar}

\begin{numar}
From  \ref{prop:incidence_cond_Cw} it follows that the Schubert cells of $\fo(E)$ are
stable under the Borel subgroup $B$ of $\spdd$ which preserves the flag
$(E_1 \subset \dots \subset E_{2n+1})$. We actually have:

\begin{prop}
The Schubert cells of $\fo(E)$ are the orbits of $B$.
\end{prop} 

\begin{proof}
The Schubert cells of $\fo(E)$ are orbits of the Borel subgroup $B_{2n+2}$
of the symplectic group $\spddd$ which preserves the flag
$(E_1 \subset \dots \subset E_{2n+1} \subset E_{2n+2} )$.
The assertion follows since $B_{2n+2}$ acts on $X_{\wz}$ trough the
restriction morphism $b \mapsto \res{b}{E}$ and the image of $B_{2n+2}$
via this morphism is $B$.
\end{proof}
\end{numar}


\begin{numar} 
We equip the subset $\wspdd \subset \wspddd$ with the order induced
by the Bruhat order of the Weyl group $\wspddd$. In this way, the 
order on $\wspdd$ describes the Bruhat order on the set of 
Schubert cells of $\fo(E)$, that is, 
for all $w,w'\in\wspdd$,
$$
C_w\subset\b{C_{w'}} \iff w\leq w'.
$$
Define the \emph{Schubert subvarieties} of the odd symplectic flag manifold
$\fo(E)$ to be the closures of the Schubert cells.
For $w \in \wspdd$, the Schubert subvariety $X_w =\b{ C_w }$ is then given by
$$X_w=\bigcup_{w'\leq w} C_{w'}.$$ 
In terms of the incidence with the flag $(E_1 \subset \dots \subset E_{2n+1})$, 
the Schubert subvariety $X_w$ is the set of flags of isotropic subspaces 
$V_{\blt}=(V_1\subset\dots\subset V_{n+1})$
which  verify
$$\dim(V_i \cap E_j) \geq r_w(i,j)$$
for all $0\leq i \leq n$ and 
$0 \leq j \leq 2n$.
\end{numar}

\begin{numar}
We can count the Schubert cells of $\foe$ of a given dimension of
and thus compute the Poincar\'e polynomial of $\foe$. We
have:
\begin{prop}  
The Poincar\'e polynomial of the variety  $\foe$ is
\begin{equation} \label{eq:poincare_fodd}   
P(\foe, q) =\frac%
        {(q^{n+1}-1) (q^{2n}-1)(q^{2n-2}-1)\dots(q^{2}-1)}%
        {(q-1)^{n+1}}.
\end{equation}
\end{prop}
Since $2, 4, \dots, 2n, n+1$ are the exponents of the Weyl group
$W(\ssD_{n+1})$ of type $\ssD_{n+1}$ (cf. \cite{Bourbaki}), it follows
that \eqref{eq:poincare_fodd}  coincides with the Poincar\'e
polynomial of the flag variety $G/B$ for $G$ of type $\ssD_{n+1}$.
We prove this proposition in \cite{eqcoh} where we also show that
$\foe$ and $G/B$ with $G$ of type $\ssD_{n+1}$ have actually the same singular
cohomology rings. 
\end{numar}
\section{The automorphism group of an odd symplectic grassmannian}

We compute here the automorphism group of the odd symplectic
grassmannian $\go(k,E)$. If $k=1$ then $\go(1,E)$ is just the
projective space $\P E$ and its automorphism group is $\PSL_{2n+1}$.
If $k=n+1$ then, as seen in
\ref{nr:orbits_in_gokdd}, the odd
symplectic grassmannian $\go(n+1,E)$ is isomorphic to the symplectic
grassmannian $\go(n,F)$, and so its automorphism group is $\pspd$.
For the remaining cases we have the following result:

\begin{prop} \label{prop:aut(gokdd)} 
For $2 \leq k \leq n$, the automorphism group of the odd symplectic
grassmannian $\gokdd$ is $\pspdd= \spdd /  \{ \pm 1  \} $.
\end{prop} 

The proof we will present here is inspired from 
a now standard 
method of computing the automorphism group of the usual grassmannian 
(see for example \cite[10.19]{Harris})
and relies on an analysis of the linear spaces contained in
 $\go(k,E)$. 

But before filling in
the details,  
let us mention that there is a different, quite straightforward way to compute
the Lie algebra of this automorphism group. This method is detailed in
\cite[\S 2.3.1]{teza} 
and roughly consists of the following steps. 
First of all, the Lie algebra of the automorphism group of
$\go(k,E)$ coincides with the space of global sections of the tangent
bundle $\hz{\go(k,E), T\go(k,E)}$  (cf. \cite[\S2.3]{Akh}).  
In order to compute this we use the description of $\goke$ as the zero
locus of $\omega $ seen as a generic section of the bundle $\Lambda ^2T^*$ on the grassmannian
$\gke$.
Then the differential of this section identifies $\Lambda ^2T^*$ with
the normal bundle of $\goke$ in $\gke$, so that the space of global
sections we want to compute is the kernel of the map
\begin{equation} \label{eq:d_omega} 
H^0\big( \goke,\res{T\gke}{\goke} \big) \pe{d\omega }{\vers} 
  H^0\big( \goke,\res{\ldt}{\goke} \big).
\end{equation}
We now use the Koszul complex associated to the section $\omega $ to
derive resolutions by locally free sheaves on the grassmannian $\gke$
for the two restrictions $\res{T\gke}{\goke}$ and 
$\res{\ldt}{\goke}$. With the help of Bott's theorem we prove that
these resolutions are acyclic, so that taking the global sections  yields
resolutions for the two spaces in \eqref{eq:d_omega}. 
Specifically, \eqref{eq:d_omega} becomes
\begin{equation}  \label{eq:d_omega_expl} 
\liesl(E) \vers \Lambda ^2E^* /\omega 
\end{equation}
where the map is sending an element $X \in \liesl(E)$ to the image in  
$\Lambda ^2E^* /\omega $ of the transform $X \cdot \omega $ 
corresponding to the $\liesl(E)$-module structure of 
$\Lambda ^2E^*$. The kernel of \eqref{eq:d_omega_expl} is then the
Lie algebra 
$\{X\in\slin(E)\mid X\cdot\omega \in \C\omega \}$
which identifies with 
$\{X\in\glin(E)\mid X\cdot\omega \in \C\omega \} / \C\Id$ and so
with 
$\liespdd=\{X\in\glin(E)\mid X\cdot\omega =0 \}$ (as Lie algebras, and
not only as vector spaces, 
since $\Id$ is central in $\liegl(E)$).
The Lie algebra of the automorphism group of $\goke$ is therefore $\liespdd$.

\begin{numar} 
We now proceed with our proof of proposition \ref{prop:aut(gokdd)}.
The outline of the proof is the following.
First we show that, in its Plücker embedding via the inclusion
$\goke \subset \gke$, all the automorphisms of $\goke$ come from
automorphisms of the ambient projective space. This implies in
particular that any automorphism of $\goke$ induces automorphisms of
the Fano schemes of linear spaces contained in $\goke$. It turns out
that a certain Fano scheme of $\goke$ is either irreducible or has
two non-isomorphic
irreducible components, and in both cases one of its irreducible
components is isomorphic to the blow-up of $\go(k-1, E)$ along its
closed $\spdd$-orbit. 
We show then that any automorphism of this Fano scheme coming from an
automorphism of $\goke$ induces further an automorphism of
$\go(k-1,E)$. Iterating this, any automorphism of $\goke$ induces an
automorphism of $\go(1,E) = \P E$, that is, an element of $\PSL_{2n+1}$,
which, as it comes from an automorphism of $\go(2,E)$, 
will be shown to be in
$\pspdd$. Finally, we show that the initial automorphism of $\goke$
coincides with the automorphism induced by this element of $\pspdd$. 
\end{numar}

\begin{numar}   
Let's start with some preliminaries. 
We first recall  Bott's theorem in the particular case
of the grassmannian $\gke$ which will be used below. 
Let $P_k \subset \sldd$ be the parabolic subgroup which
preserves the space generated by the last $k$ vectors of the basis
$ \{ e_0, \dots, e_{2n}  \} $, so that $\gke= \sldd / P_k$.  In
$\sldd$ we fix the maximal torus of diagonal matrices in the basis
$ \{ e_0, \dots, e_{2n}  \} $ and the Borel subgroup of upper
triangular matrices in the same basis; the parabolic $P_k$ contains
then the opposite Borel subgroup. 
The weights of $\sldd$ identify with the tuples 
$\lambda =(\lambda _0,\dots,\lambda _{2n})$ 
modulo shifts
$(\lambda _0,\dots,\lambda _{2n})\mapsto
(\lambda _0+a,\dots,\lambda _{2n}+a)$.
A weight $\lambda =(\lambda _0,\dots,\lambda _{2n})$ 
is \emph{singular} if it has two equal parts and
\emph{dominant} if 
$\lambda_0\geq\cdots\geq\lambda _{2n}$. 
In the latter case, if the $\lambda _i$ are non negative, we see
$\lambda $ as a partition and denote $|\lambda |=\sum\lambda _i$ its
\emph{weight} and $\ell(\lambda )$, the number of non zero parts, its \emph{length}.
For $\lambda $ dominant, 
the simple $\sldd$-module with highest weight $\lambda $ is isomorphic
to the Schur power $S_\lambda E$.

Irreducible homogeneous vector bundles on $\gke$ are determined by
irreducible $P_k$-modules, which in turn are determined by 
$P_k$-dominant weights, that is weights 
$\lambda =(\lambda _0, \dots,\lambda _{2n})$ 
such that 
$\lambda_0 \geq\cdots\geq \lambda _{2n-k}$ and
$\lambda _{2n-k+1} \geq\cdots\geq \lambda _{2n}   $. 
If $T$ and $Q$ denote the tautological, respectively the quotient,
bundles on the grassmannian $\gke$, whose fibers at a point 
$V \in \gke$ 
are $V$, respectively $E/V$, then the irreducible bundle corresponding
to a $P_k$-dominant weight 
$\lambda =(\lambda _0,\dots,\lambda _{2n})$ 
is given in terms of Schur powers of $T$ and $Q$ by 
$$
E_{\lambda }=S_{\lambda '}Q\otimes S_{\lambda ''}T,
$$
where  
$\lambda '=(\lambda_0\geq\cdots\geq\lambda _{2n-k})$ and
$\lambda ''=(\lambda _{2n-k+1}\geq\cdots\geq\lambda _{2n})$.

The Weyl group of $\sldd$ is isomorphic to the permutation group
$S_{2n+1}$ and acts on weights by 
$w \cdot (\lambda _0,\dots,\lambda _{2n}) = (\lambda _{w^{-1} (0)},\dots,\lambda _{w^{-1} (2n)}).$  
Let $\rho =(2n+1,2n,\dots,1)$ denote the smallest strictly decreasing
partition, 
which modulo a shift coincides with the half-sum of the
positive roots.
For $\lambda $ a $P_k$-dominant weight, the cohomology of the
associated irreducible bundle $E_\lambda $ is described by Bott's theorem:
\begin{thm*}
If $\lambda +\rho $ is singular then $H^i(\gke,E_{\lambda})=0$ for all $i$.
If not, let  $w\in S_{2n+1}$ be the unique permutation such that
$w(\lambda +\rho )$ is dominant,  and let $\ell$ be its length.
Then $H^{\ell} (\gke,E_{\lambda})$ is the
simple $\sldd$-module  
of highest weight $w(\lambda +\rho )-\rho $ and 
$H^i(\gke,E_{\lambda})=0$ for all $i\neq \ell$.
\end{thm*}

We'll also need the plethysm formula for the decomposition of 
$\Lambda^j(\Lambda ^2 V)$ as a $\GL(V)$-module. 
Recall first that the Frobenius notation 
$\lambda =(a_1,\dots,a_r \mid b_1,\dots,b_r)$
for a partition $\lambda $ encodes its \emph{hook decomposition}, that
is, in
the Ferrers diagram of $\lambda $ there are $a_i$ boxes to the right
of the $i$-th box on the diagonal and $b_i$ boxes below. The
\emph{rank} $r$ of $\lambda $ is the number of boxes on the diagonal.
We have (cf. \cite{Macdonald}) 
\begin{equation} \label{eq:pleth} 
\Lambda ^i(\Lambda ^2V)=\bigoplus_{|\lambda |=i}S_{\lambda ^-}V
\end{equation}
where the sum is indexed by all strictly decreasing partitions of
weight $i$ and for a strictly decreasing partition $\lambda $ we
denote
$$
\lambda ^-=(\lambda _1-1,\dots,\lambda _{\ell}-1 \mid \lambda _1,\dots,\lambda _{\ell}).
$$
\end{numar}

\begin{numar} We embed $\goke$ in $\gke$ and further, via the Plücker
embedding of $\gke$, in $\P \Lambda ^kE$. We set out to prove that
any automorphism of $\goke$ comes from an automorphism of the ambient
projective space $\P \Lambda ^kE$.

\begin{prop} \label{prop:lin_normale} 
The embedding $\goke \subset \P\Lambda ^k E$ is linearly
normal, that is, 
the restriction morphism
$$
H^0\big(\P \Lambda ^k E,\cO_{\P\Lambda ^k E}(1)\big)\vers
H^0\big(G_{\omega }(k,E),\cO_{G_{\omega }(k,E) }(1)\big)
$$
is surjective.
\end{prop}

\begin{proof}
To simplify the notation, we write $G$ instead of $\gke$ and $\go$
instead of $\goke$. The line bundle $\cO(1)$ on $G$ corresponding to
the Plücker embedding is isomorphic to $\det T^*$. Since 
$H^0(G, \det T^*)=\Lambda ^k E^*$, it suffices to show that the
restriction morphism 
\begin{equation}  \label{eq:restriction_G_Go} 
H^0(G, \cO_G(1)) \vers H^0(\go, \cO_{\go} (1))
\end{equation}
is surjective.
We regard $\go$ as the zero locus of $\omega $ seen as a section of
the bundle $\Lambda ^2 T^*$ on $G$. We consider the Koszul complex of
$\omega $ which we twist by $\det T^*$ and obtain the exact complex
\begin{multline} \label{eq:koszul_detT*}
0\vers \Lambda^r(\Lambda^2T)\otimes \det T^* \vers \cdots \vers \Lambda ^2(\Lambda ^2T)\otimes \det T^* \vers \\
\vers \Lambda ^2 T\otimes \det T^* \vers \det T^*\vers \res{\det T^*}{\go}\vers 0
\end{multline}
where $r=\frac12 k(k-1)$. To obtain the surjectivity of 
\eqref{eq:restriction_G_Go} 
it suffices to show that this complex is acyclic, ie that the bundles
in \eqref{eq:koszul_detT*} have no higher cohomology (this follows
in a standard way by chopping \eqref{eq:koszul_detT*} into short exact
sequences and using the corresponding cohomology long exact sequences).
By \eqref{eq:pleth}, we have
$$ \Lambda ^j(\Lambda ^2T)\otimes \det T^* = \bigoplus_{|\lambda |=j}S_{\lambda  ^-}T\otimes \det T^*$$
so any irreducible component of
$\Lambda ^j(\Lambda ^2T)\otimes \det T^*$ 
is of the form
$E_\eta=S_\mu T$ where $\mu =(\mu _1\geq\cdots\geq\mu _k)$ 
is the skew partition 
$(\lambda ^-_1-1, \dots, \lambda ^- _k-1)$
for some strictly decreasing partition $\lambda $ of weight $j$.
Then
$\eta=(0,0,\dots,0,\mu _1,\dots,\mu _k)$
and
$$
\eta+\rho=(2n+1,2n,\dots,k+1,\mu _1+k,\dots,\mu _k+1).
$$
Consider first the case $j \geq 2$. By construction, we have
$\ell(\lambda ^-)=\lambda _1+1$ and $\lambda ^-_1=\lambda _1$.
As $T$ is of rank $k$, $\ell(\lambda ^-) \leq k$, 
since otherwise $S_{\lambda ^-}T=0$, 
so $\lambda _1\leq k-1$.
Therefore
$\mu _1+k  = \lambda _1 - 1 + k \leq 2k-2 \leq 2n+1.$
On the other hand, since $| \lambda |=j \geq 2$ and $\lambda $ is
strictly decreasing, we have $\lambda _1 \geq 2$ and so
$ \mu _1+k  = \lambda _1 - 1 + k \geq k+1$.
The partition $\eta+\rho$ has therefore two equal parts and so is
singular. By Bott's theorem, it follows that
$$
H^i(G,\Lambda ^j(\Lambda ^2T)\otimes \det T^*)=0\text{ for all }i.
$$
For $j=1$, the bundle
$\Lambda ^2 T\otimes \det T^* =\Lambda ^{k-2}  T^*$ is given by a
dominant weight and therefore has no higher cohomology.
Similarly for $j=0$ and the bundle $\det T^*$, and we are done.
\end{proof}

\begin{prop} \label{prop:aut_ambient} 
Any  automorphism of the odd symplectic grassmannian $\goke$ is induced
by an automorphism of the ambient projective space
$\P\Lambda ^kE$.
\end{prop}

\begin{proof} 
This follows by a standard argument 
(as in \cite[Example 7.1.1]{Hartshorne}, for example)
from the fact that $\goke$ is linearly normal in $\P\Lambda ^kE$ and
its Picard group is $\Z$.
\end{proof}
\end{numar}

\begin{numar} 
Denote by $\Lambda ^{\Spc{k}}E$ the kernel of the contraction
$\llcorner\omega : \Lambda ^kE \vers \Lambda ^{k-2}E$.
Since 
$H^0(G,\Lambda ^{k-2} T^*) = \Lambda ^{k-2}E^* $, taking global
sections in \eqref{eq:koszul_detT*} and using the vanishings from the
proof of  \ref{prop:lin_normale}, we get the exact sequence
$$
0\vers \Lambda ^{k-2}E^*  \pe{\wedge \omega }{\vers}    \Lambda^kE^*  
     \vers  H^0(\go,\cO_{\go}(1))  \vers 0.
$$
Taking the duals, we obtain 
$H^0(\go,\cO_{\go}(1))^* = \Lambda ^{\Spc{k}}E$ so $\goke$ embeds in
the projective subspace 
$\P\Lambda ^{\Spc{k}}  E \subset   \P\Lambda ^k E$ as a non-degenerate
subvariety.
We actually have:

\begin{prop}  \label{prop:sect_lin}  
$G_{\omega} (k,E)=G(k,E)\cap\P\Lambda ^{\Spc{k} }E.$
\end{prop} 

\begin{proof} 
Let  $V\in G(k,E)$ and $\{v_1,\dots,v_k\}$ a basis of  $V$. The contraction of 
$v_1\wedge \dots\wedge v_k$
by $\omega $ is explicitly given by
$$
(v_1\wedge \dots\wedge v_k)  \llcorner  \omega 
    =\sum_{1\leq i<j\leq k}(-1)^{i+j}\omega (v_i,v_j)
v_1 \wedge\dots\wedge  \widehat{v_i}
            \wedge\dots\wedge
                 \widehat{v_j} \wedge\dots\wedge  v_k
$$  
so
$v_1\wedge \dots\wedge v_k$ is in the kernel of  $\llcorner\omega $ if
and only if $V$ is isotropic.
\end{proof}

Note, though, that the intersection in \ref{prop:sect_lin}  is not transverse.
Let us denote by $\Aut(\P\Lambda ^{\Spc{k}}E,G_\omega (k,E))$ the group of
automorphisms of the projective space  $\P\Lambda ^{\Spc{k}}E$
which preserve  
$G_\omega (k,E)$. 
We have natural morphisms
$$
\pspdd \vers \Aut\big(\P\Lambda ^{\Spc{k}}E,G_\omega (k,E)\big)
   \vers  \Aut\big(G_\omega (k,E)\big).
$$ 
The second morphism is an isomorphism, since it is surjective by
\ref{prop:aut_ambient} and injective since $\goke$ is non-degenerate
in $\P\Lambda ^{\Spc{k}}E$.
The first morphism is injective, since any line in $E$ is recovered as
the intersection of all the $k$-dimensional isotropic subspaces which
contain it.
Therefore, to conclude our proof of  \ref{prop:aut(gokdd)} it suffices
to show that the first morphism is surjective.

We start analyzing the linear spaces contained in $\goke$. We recall
first the description we have in the case of the grassmannian  $\gke$
(cf.  \cite[6.9]{Harris}):

\begin{prop} \label{prop:espaces_gkdd}  
A maximal linear space contained in the grassmannian $\gke$ is either
of the form
$$
\{V\in G(k,E)\mid V \supset V^{k-1}\}
$$
for some subspace $V^{k-1}\subset E$
of dimension $k-1$, or of the form
$$
\{V\in G(k,E)\mid V \subset V^{k+1}\}
$$ 
for some subspace $V^{k+1}\subset E$
of dimension $k+1$. 
\end{prop}

For the odd symplectic grassmannian $\goke$ we have:

\begin{prop} \label{prop:lin_max} 
The maximal linear spaces contained in $\goke$ are of one of the
following two types:
\begin{align*} 
&\text{type I: } \quad\big\{V \in \goke \mid V^{k-1}\subset V\subset (V^{k-1})^{\bot}\big\}=%
  \begin{cases}
         \P^{2(n-k)+2}&\text{if }V^{k-1}\not\ni e_0\\
         \P^{2(n-k)+3}&\text{if }V^{k-1}\ni e_0
  \end{cases}
\\
\intertext{for an isotropic subspace $V^{k-1}\subset E$ of dimension $k-1$,} 
&\text{type II: } \quad\big\{V\in \goke \mid V\subset V^{k+1}\big\}=%
  \begin{cases}
        \P^{k}&\text{if }V^{k+1}\text{ isotropic}\\
        \P^{1}&\text{if } \rank(\res{\omega }{V^{k+1}})=2
  \end{cases}
\end{align*}
for a subspace $V^{k+1}\subset E$ of dimension $k+1$.
\end{prop} 

\begin{proof} 
Since the odd symplectic grassmannian $\goke$ is cut in $\gke$ by a
linear space (cf.  \ref{prop:sect_lin}),
the maximal linear spaces contained in $\goke$ are simply the
intersections of $\goke$ with the maximal linear spaces contained in $\gke$.

Let $L$ be a maximal linear space contained in  $\gke$ of the first
kind described in  \ref{prop:espaces_gkdd}, that is,
$L=\{V\in G(k,E)\mid V \supset V^{k-1}\}$ for some subspace $V^{k-1}\subset E$
of dimension $k-1$. 
If  $V^{k-1}$ is not isotropic then $L\cap\goke=\emptyset$. 
Suppose that $V^{k-1}$ is isotropic and let  $V$ be a point of $L$.
If $V$ is  isotropic then $V\subset (V^{k-1})^{\bot}$ and, conversely,
if  $V\subset (V^{k-1})^{\bot}$ then $V$ is isotropic, so
$$
L\cap\goke=\big\{V\in G(k,E)
\mid V^{k-1}\subset V\subset (V^{k-1})^{\bot}\big\}\simeq 
\P\big((V^{k-1})^{\bot}/V^{k-1}\big).
$$ 
If $V^{k-1}$ does not contain the kernel $R$ of $\omega $ then
$\dim  (V^{k-1})^{\bot} = 2n-k+2$
and  $\P((V^{k-1})^{\bot}/V^{k-1})   \simeq  \P^{2n-2k+2}$.
If $V^{k-1}$ contains $R$ then
$\dim  (V^{k-1})^{\bot}= 2n-k+3$
and  $\P((V^{k-1})^{\bot}/V^{k-1})  \simeq   \P^{2n-2k+3}$. 

Let now $L$ be a maximal linear space contained in $\gke$ of the
second kind described in \ref{prop:espaces_gkdd}, that is,
$L=\{V\in G(k,E)\mid V\subset V^{k+1} \}$
for some subspace
$V^{k+1}\subset E$ 
of dimension $k+1$. 
The intersection $L\cap\goke$ is then
the grassmannian of isotropic hyperplanes in $V^{k+1}$.
By the lemma below
we have three possibilities, according to the rank of 
$\res{\omega }{V^{k+1}}$: 
If $V^{k+1}$ is isotropic then
$$
L\cap\goke=L=\{V\in G(k,E)\mid V\subset V^{k+1} \}\simeq\P^k;
$$
If $\rank( \res{\omega }{V^{k+1}})=2$ then
$$
L\cap\goke=\big\{V\in G(k,E)\mid \Ker(\res{\omega }{V^{k+1}})\subset
V\subset V^{k+1}\big\}\simeq\P^1;
$$
If $\rank( \res{\omega }{V^{k+1}})>2$ then  $L\cap\goke=\emptyset$.
This concludes the proof.
\end{proof}

\begin{lema} 
Let $W$ be a space equipped with a non zero skew-form $\omega $. 
If  $\rank \omega =2$ then the isotropic  hyperplanes of $W$ are
precisely those which contain $\Ker \omega $. If $\rank \omega >2$ then there is no isotropic hyperplane in $W$.
\end{lema}

\begin{proof}
Denote still by
$\omega $ the canonical morphism
$\omega :W\to W^*$. 
We write  $\bot_\omega$ to denote the orthogonal with respect to $\omega $ 
and $\bot$ to denote the orthogonal with respect to the duality $W\otimes W^*\to\C$.
Let  $H\subset W$ be a hyperplane. Then
$$
H\text{ isotropic}\iff H\subset H^{\bot_\omega}=\omega (H)^{\bot}\iff
H^{\bot}\supset \omega (H).
$$    
Since  $\dim H^{\bot}=1$, this is equivalent to $\omega (H)=0$ or $\omega (H)=H^{\bot}.$
We cannot have $\omega (H)=0$ since $\Ker\omega $ is of codimension at
least  2.
If $\omega (H)=H^{\bot}$   then  $\Ker\omega \cap H$ is a  hyperplane of $H$ so 
$\codim\Ker\omega \leq2$, and so $\codim\Ker\omega =2$ and $\Ker\omega \subset H$.
\end{proof}

\end{numar}

\begin{numar} \label{nr:F(gokdd)} 
The linear spaces contained in $\goke$ which will prove useful in our
computation of the automorphism group of $\goke$ are those of
dimension $2(n-k)+2$. 
Recall that for a projective variety $X\subset \P^N$ and an integer
$\ell\in\N$, the \emph{Fano scheme}
$\bF(\ell,X)$ is the subvariety of the grassmannian
$\bbG(\ell,N)$ of  $\P^{\ell}$'s in $\P^N$ parametrizing
the linear subspaces of $\P^N$ of dimension $\ell$ contained in $X$.

In what follows we will focus on the Fano scheme
$\bF(2(n-k)+2,\goke)$
of linear spaces of dimension 
$2(n-k)+2$ contained in $\goke$. For simplicity, we will denote
it $\bF$.

\begin{prop} \label{prop:sch_fano}
If $k < \frac23(n+1)$  then the Fano scheme $\bF$ parametrizing the 
$\P^{2(n-k)+2}$'s contained in $\goke$ is irreducible and isomorphic
to the
variety
$$
\bF^1=\big\{(V^{k-1},H)\in G_\omega (k-1,E)\times G(2n-k+2,E)\;\big|\;
               V^{k-1}\subset H\subset(V^{k-1})^{\bot}\big\}
$$
which identifies to the blow-up of 
$\go(k-1,E)$ along its closed $\spdd$-orbit.

If $k \geq \frac23(n+1)$ then $\bF$ has two disjoint irreducible components: 
one is isomorphic to the variety $\bF^1$ above, 
the other is isomorphic to the variety
$$
\bF^2=\big\{(W,V^{k+1})\in G(3k-2(n+1),E)\times G_\omega (k+1,E)
               \;\big|\;   W\subset V^{k+1}\big\}
$$
and these two components are not isomorphic.
\end{prop} 

We will prove this proposition in several steps. First, a little piece
of terminology. 
Any $\P^{2(n-k)+2}$ contained in $\goke$ is included in some maximal
linear space contained in $\goke$.
If the maximal linear space is of type I (according to
\ref{prop:lin_max})  then we will say that the
$\P^{2(n-k)+2}$ is of the \emph{first kind}. Otherwise, ie if the
maximal linear space is of type II,
we will say
the $\P^{2(n-k)+2}$ is of the \emph{second kind}.
Note that this terminology is coherent, that is, a 
$\P^{2(n-k)+2}$ cannot simultaneously be
 of the first and of the
second kind since the intersection of a type I maximal linear space
with a type II maximal linear space is 
contained in a $\P^1$:
\begin{multline*} 
\big\{V\in \goke \mid V^{k-1}\subset V\subset (V^{k-1})^{\bot}\big\} \cap 
\big\{V\in \goke \mid V\subset V^{k+1}\big\} \subset \\
\subset \big\{V\in G(k,E)\mid V^{k-1}\subset V\subset V^{k+1}\big\} \simeq \P^1.
\end{multline*}
Note also that if
$k <  2(n-k)+2$, or equivalently $k < \frac23 (n+1)$, then from the
definition it follows that there is no
$\P^{2(n-k)+2}$ of the second kind.

\begin{lema} \label{lema:F1} 
The $\P^{2(n-k)+2}$'s of the first kind are parametrized by the
variety $\bF^1$.
\end{lema}

\begin{proof} 
Consider a
$\P^{2(n-k)+2}$ 
of the first kind. If it is maximal, then according to 
\ref{prop:lin_max} it is determined by giving 
a $V^{k-1}$ isotropic of dimension $k-1$ which does not contain $e_0$.
If it is not maximal then it is a hyperplane in a maximal
$\P^{2(n-k)+3}$ of  type I, ie it is determined by giving a
$V^{k-1}$ isotropic of dimension $k-1$ which contains $e_0$ and a hyperplane $H$ of 
$(V^{k-1})^{\bot}$ which contains $V^{k-1}$ :
$$
\P^{2(n-k)+2}=\big\{V\in G(k,E)\;\big|\; V^{k-1}\subset V\subset
                                               H\subset (V^{k-1})^{\bot}\big\}.
$$
Now consider a point $(V^{k-1},H) \in \bF^1$. 
If  $V^{k-1}\not\ni e_0$ then $\dim (V^{k-1})^{\bot} =2n-k+2$, so
$H=(V^{k-1})^{\bot}$.
If $V^{k-1}\ni e_0$ then $\dim (V^{k-1})^{\bot} =2n-k+3$ 
and so $H$ is a hyperplane in $(V^{k-1})^{\bot}$.
This shows that any $\P^{2(n-k)+2}$ 
of the first kind is of the form
\begin{equation} \label{eq:first_kind} 
 \{V\in G(k,E)\mid V^{k-1}\subset V\subset H\}
\end{equation}
for some  $(V^{k-1},H) \in \bF^1$. 
Conversely, any linear space of the form \eqref{eq:first_kind} is a 
$\P^{2(n-k)+2}$ contained in $\goke$, and we are done.
\end{proof}

\begin{lema}
The  $\P^{2(n-k)+2}$'s of the second kind are parametrized by 
$\bF^2$.
\end{lema}

\begin{proof} 
A $\P^{2(n-k)+2}$ of the second kind is included in a maximal
$\P^k$ of type II, so, according to \ref{prop:lin_max}, 
it is defined by giving a $V^{k+1}$ isotropic
dimension $k+1$ 
and a subspace $V^{3k-2(n+1)}\subset V^{k+1}$ of dimension $3k-2(n+1)$:
\begin{equation}  \label{eq:second_kind}  
\P^{2(n-k)+2}=\{V\in\gke  \mid V^{3k-2(n+1)}\subset V\subset V^{k+1} \}.
\end{equation}
Conversely, any linear space of the form \eqref{eq:second_kind} is 
contained in $\goke$.
\end{proof}

\begin{lema} \label{prop:F1=eclatement} 
The projective variety $\bF^1$ is  isomorphic to the blow-up of the
odd symplectic grassmannian
$\go(k-1, E)$ along its closed  $\spdd$-orbit.
\end{lema}

\begin{proof} 
Let $p:\bF^1\to G_\omega (k-1,E)$ denote the projection on the first
factor, and let $V^{k-1} \in  \go (k-1,E)$.
From the  discussion in the proof of lemma \ref{lema:F1}, it follows that if
$V^{k-1}\not\ni e_0$ then the fiber $p^{-1} ( V^{k-1} )$ 
reduces to 
the point $( V^{k-1}, (V^{k-1})^{\bot} )$  
and if $V^{k-1}\ni e_0$
then the fiber $p^{-1} ( V^{k-1} )$ identifies with
$$
\{  H  \mid H\text{ hyperplane in }(V^{k-1})^{\bot}, H\supset V^{k-1} \}
           =\P((V^{k-1})^{\bot}/V^{k-1})^* .
$$
So, above the open $\spdd$-orbit 
$X_1\subset \go(k-1,E)$
the projection $p$ is an isomorphism, and above the closed 
$\spdd$-orbit $X_0 \subset \go(k-1,E)$
it is a fibration with fiber
$\P^{2n-2k+3}$.

Since $X_0$ is of codimension $2n-2k+4$ in $\go (k-1,E)$ it follows
that $p^{-1} (X_0)$ is a divisor in $\bF^1$, and to prove that
$\bF^1$ is the blow-up of $X_0$ in $\go (k-1,E)$ it suffices to show
that $\bF^1$ is smooth  (cf. \cite[Ch. 4,\S 6]{GH}).
We will show that 
$\bF^1$ is the zero locus of a sufficiently generic section of a
vector bundle on the partial flag variety
$$
\bbF(k-1,2n-k+2,E) =
 \big\{(V^{k-1},V^{2n-k+2})\mid V^{k-1}\subset V^{2n-k+2}\subset E\big\}
$$
(with the convention, implicitly used already, to write dimensions as
superscripts).
Note that this flag variety is a component of the Fano scheme
$\bF(2(n-k)+2,G(k,E))$ parametrizing the $\P^{2(n-k)+2}$'s contained
in the grassmannian $\gke$.

For simplicity we denote in what follows $a=k-1$ and $b=2n-k+2$. 
Let  $T^a$ and $T^b$ be the  tautological bundles on the flag variety
$\bbF(a,b,E)$ 
whose fibers at a point
$V^\blt=(V^a\subset V^b) \in \bbF(a,b,E)$
are $V^a$, respectively $V^b$. 
The skew-form $\omega $ defines a section of the bundle
$\tast\otimes\tbst$ whose value at the point $(V^a\subset V^b)$ is the
restriction of 
$\omega $ to the subspace $V^a\times V^b \subset E \times E$.

A point $(V^a\subset V^b)$ of $\bbF(a,b,E)$ is in $\bF^1$ if and only
if  $V^a$ is isotropic and $V^b\subset(V^a)^{\bot}$, which is
equivalent to $V^b\subset(V^a)^{\bot}$
since $V^a\subset V^b$.
In other words, $\bF^1$ is the zero set of the section
$\omega \in H^0(\bbF(a,b,E),\tast\otimes\tbst)$.
But the bundle $\tast\otimes\tbst$ is not suitable for our purpose
since its rank is too big. We replace it by a subbundle
of which $\omega $ is still a section.
Let $\E$ be the kernel of the surjective morphism  
$\tast\otimes\tbst \to S^2\tast$ given by the
composition $\tast\otimes\tbst \to \tast\otimes\tast \to S^2\tast$.
Then sections of $\tast\otimes\tbst$ which come from $\Lambda ^2E^*$
are sections of $\E$, as is the case in particular for $\omega \in \Lambda ^2E^*$.
The fiber of $\E$ at a point $(V^a\subset V^b)$ identifies, in a
basis adapted to the flag $V^a\subset V^b$, 
to $a\times b$ matrices of the form
$$
\begin{pmatrix} A & B 
\end{pmatrix}, \quad A\text{ antisymmetric.}
$$  
Such a matrix can always be completed to a 
$(2n+1)\times(2n+1)$
antisymmetric matrix
$$
\begin{pmatrix}
A & B & * \\
* & * & * \\
* & * & *
\end{pmatrix} 
$$
which shows that $\E$ is generated by its global sections which come
from $\Lambda ^2E^*$.
Since $\omega $ is generic in $\Lambda ^2E^*$ and $\bF^1$ is non
empty, we conclude by Bertini's theorem.
\end{proof}

The projective variety $\bF^2$ is a fibration in grassmannians
$G(3k-2(n+1),k+1)$ over the odd symplectic grassmannian
$\go(k+1,E)$, it is then smooth, irreducible, of dimension
\begin{equation*} 
\begin{split} 
\dim\bF^2 &=\dim G_\omega (k+1,E) +\dim G(3k-2(n+1),k+1) \\
          &=(k+1)(2n-k)-\binom{k+1}{2}+[3k-2(n+1)][2(n-k)+3].
\end{split}
\end{equation*}
On the other hand, the variety $\bF^1$ is birational to the odd
symplectic grassmannian $\go (k-1,E)$, so it is irreducible,
of dimension 
\begin{equation*} 
\dim\bF^1 =\dim G_\omega (k-1,E)=(k-1)(2n-k+2)-\binom{k-1}{2}.  
\end{equation*}
We have then
\begin{equation*} 
\dim\bF^2- \dim\bF^1 = [3k-2(n+1)][2(n-k)+1]-1
\end{equation*}                   
so $\dim \bF^1 = \dim \bF^2$ only when 
$n=k=3$, in which case  $\bF^1$ and $\bF^2$ 
are both of dimension 9.

\begin{lema}
When $n=k=3$, the varieties $\bF^1$ and $\bF^2$ are not isomorphic.
\end{lema}

\begin{proof} 
We have
\begin{gather*}  
\bF^1=\big\{(V^2,H^5)\in G_\omega (2,7)\times G(5,7) \;\big|\;
    V^2\subset H\subset (V^2)^\bot \big\} \\
\bF^2=\big\{(W^1,V^4)\in G(1,7)\times G_\omega (4,7) \;\big|\;
    W^1\subset V^4 \big\}.
\end{gather*}
The variety $\bF^1$ is the blow-up of  $\go(2,7)$ along the closed
$\Sp_7$-orbit which is isomorphic to $\go (1,6)\simeq \P^5$, so
the exceptional divisor of $\bF^1$ is a fibration
in $\P^3$ over $\P^5$.
The variety $\bF^2$ is a fibration in $G(1,4)\simeq\P^3$ over the odd
symplectic grassmannian $\go(4,7)$.
Note that the varieties $\bF^1$ and $\bF^2$ are both smooth of
dimension 9 and of Picard group $\Z^2$. We show they don't have the
same singular cohomology.

If $p:\bF^1\to \go(2,7)$ is the canonical projection and
$D=p^{-1}(\P^5) $ the exceptional divisor,
then, as $\Z$-modules, we have (cf. \cite[Ch. 4,\S 6]{GH})
$$
H^*(\bF^1,\Z)\simeq H^*(\go(2,7),\Z)\;\oplus\;H^*(D,\Z)/H^*(\P^5,\Z).
$$
The $\Z$-module $H^*(\go(2,7),\Z)$ is free, generated by the classes
of the Schubert subvarieties of $\go (2,7)$. It is of rank 
$\#\indiceadm{2}{7}=18$.
Since $p:D\to\P^5$ is a fibration in $\P^3$ over $\P^5$, the ring
$H^*(D,\Z)$ is a free  $H^*(\P^5,\Z)$-module of rank 4
(cf. \cite{GH}). As  $H^*(\P^5,\Z)$ is a free $\Z$-module of rank 6,
we obtain that
$H^*(D,\Z)/H^*(\P^5,\Z)$ is a free $\Z$-module of rank 18, and so
$H^*(\bF^1,\Z)$ is of rank 36.

Similarly, $\bF^2$ is a fibration in $\P^3$ over $\go(4,7)$ so
$H^*(\bF^2,\Z)$ is a free module of rank 4 over $H^*(\go(4,7),\Z)$.
The latter is a free $\Z$-module of rank 
$\#\indiceadm{4}{7}=8$ so $H^*(\bF^2,\Z)$ is of rank 32, and so
$\bF^1$ and $\bF^2$ are not isomorphic.
\end{proof}

This concludes
the proof of the proposition \ref{prop:sch_fano}. 
\end{numar}

\begin{numar} 
We can tackle now the final step in our calculation of the
automorphism group of $\goke$. Recall that what we are left to prove
now is that any element of
$\Aut(\P\Lambda ^{\Spc{k}}E,G_\omega (k,E))$
comes from $\pspdd$.

An element
$g \in \Aut(\P\Lambda ^{\Spc{k}}E,G_\omega (k,E))$
induces an automorphism of the Fano scheme $\bF$ of \ref{nr:F(gokdd)}.
By \ref{prop:sch_fano},  $\bF$ is either irreducible and isomorphic to $\bF^1$, 
or it has a second irreducible component  $\bF^2$ which is not
isomorphic to $\bF^1$. Therefore $g$ induces an 
automorphism of $\bF^1$, which is to say that $g$ sends 
$\P^{2(n-k)+2}$'s of the first kind to $\P^{2(n-k)+2}$'s of the first kind.

Let $p:\bF^1\to \go(k-1,E)$ be the projection.

\begin{prop}
The automorphism induced by $g$ on $\bF^1$ permutes the fibers of $p$
and induces an automorphism $\tilde g$ of $\go(k-1,E)$.
\end{prop}

\begin{proof} 
We abuse notation and denote the automorphism induced by $g$ on
$\bF^1$ by the same letter $g$.
The projection $p$ is an isomorphism over the open $\spdd$-orbit $X_1$
of $\go(k-1,E)$ (cf. the proof of lemma \ref{prop:F1=eclatement}) and
the points of $p^{-1}(X_1)$ correspond to the $\P^{2(n-k)+2}$'s of the
first kind which are maximal as linear spaces contained in $\goke$.
As $g$ sends a maximal $\P^{2(n-k)+2}$ of the first kind on a maximal
$\P^{2(n-k)+2}$ of the first kind, it follows that $g$ sends  the open
set $p^{-1}(X_1)$ on itself.

Let now $(V_1^{k-1} ,H_1)$, $(V_2^{k-1} ,H_2)$ be two points of 
$p^{-1}(X_0)\subset \bF^1$ and $L_1$, $L_2$ the corresponding
$\P^{2(n-k)+2}$'s, that is, 
\begin{align*}
L_1 &=\{V\in G(k,E)\mid V_1^{k-1}\subset V\subset H_1\}  \\  
L_2 &=\{V\in G(k,E)\mid V_2^{k-1}\subset V\subset H_2\}.   
\end{align*}
If $V_1^{k-1}=V_2^{k-1}=V^{k-1}$, then
$H_1\cap H_2$ is of codimension 2 in $(V^{k-1}   )^{\bot}$ and $L_1$ and $L_2$ 
cut along a $\P^{2(n-k)+1}$:
$$
L_1\cap L_2
   =\{V\in G(k,E)\mid V^{k-1}\subset V\subset H_1\cap H_2\}
     \simeq  \P^{2(n-k)+1}.
$$
If $V_1^{k-1}\neq V_2^{k-1}$ then the intersection
$$
L_1\cap L_2
   =\{V\in G(k,E)\mid 
         V_1^{k-1}+V_2^{k-1}\subset V\subset H_1\cap H_2\}
$$
is either empty or it reduces to a point: if $V\in L_1\cap L_2$ then necessarily
$V=V_1^{k-1}+V_2^{k-1}$ since  $V_1^{k-1}$ and $V_2^{k-1}$ are
distinct hyperplanes in $V$.
It follows that $g$ sends a fiber of $p$ on a fiber of $p$, and so it
induces a bijection $\tilde g$ such that the following square 
commutes:
\begin{equation} \label{eq:carre_g_tilde}
\begin{array}{c}  
\xymatrix{%
\bF^1 \ar[r]^{g} \ar[d]_{p} & \bF^1  \ar[d]^{p} \\
\go(k-1,E) \ar[r]^{\tilde g}  & \go(k-1,E).
}
\end{array}
\end{equation}
Since $\bF^1$ is projective, the projection $p$ is a closed map and so
$\tilde g$ is continuous for the Zariski topology.
On the other hand, $\tilde g$ is regular on the open orbit  $X_1$
since $p$ is an isomorphism over $X_1$.
The closed orbit $X_0$ is of  codimension at least  2, 
and since $\go(k-1,E)$ is normal this implies that 
$\tilde g$ is regular everywhere. Indeed, for any affine open set 
$U \subset  \go(k-1,E)$ the preimage $\tilde g ^{-1} (U)$ is normal
and the restriction $\res{\tilde g}{ \tilde g ^{-1} (U) }$ is given
by some rational functions, which are regular
outside a codimension 2 closed set and
therefore extend to the whole of $\tilde g ^{-1} (U)$. It follows that
$\res{\tilde g}{X_1}$ extends to a regular map on $\go(k-1, E)$ which,
by continuity, has to be $\tilde g$.
\end{proof}

We have just shown that an automorphism of $\goke$ induces an
automorphism of $\go(k-1,E)$. Iterating, it follows that an
automorphism of $\goke$ induces an automorphism of
$\go(1,E) = \P E$, that is, an element  $h \in \PSL(E)$. 
We now only need to show that $h$ is in $\pspdd$ and that the
automorphism of $\goke$ it defines coincides with our initial automorphism.
This is done in the following two lemmas.

\begin{lema}
If the automorphism of $G(2,E)$ induced by some $h\in\PSL(E)$ preserves $\go(2,E)$
then $h \in \pspdd$.
\end{lema}  

\begin{proof} 
Let $\hat h$ be a representative of $h$ in $\SL(E)$. Then $\hat h$
preserves orthogonality, that is, for all $u,v \in E$
$$
u\bot v \Rightarrow \hat hu\bot \hat hv,
$$
since an isotropic subspace of dimension 2 is generated by two
orthogonal vectors and conversely. In other words, the automorphism
$(u,v)\mapsto (\hat hu,\hat hv)$ preserves the affine closed set
$$
\{(u,v)\in E\times E\mid \omega (u,v)=0\}.
$$
Since 
$\omega $ is irreducible as a degree two polynomial on $E \times E$,
it follows that there is $\lambda \in\C^*$ such that
$\omega (\hat hu, \hat hv)=\lambda  \omega (u,v)$
for all $u,v\in E$, ie that
$\hat h\omega =\lambda^{-1} \omega$.
So, if  $\mu \in \C^*$ is a square root of $\lambda^{-1}$, then
$(\mu \hat h)\omega =\omega $, that is,
$\mu \hat h\in \spdd$ and so $h \in \pspdd$.
\end{proof}

\begin{lema}
Let  $g$ be an automorphism of $\goke$ and $\tilde g$ the automorphism
it  induces on $\go(k-1,E)$. If $\tilde g$ is induced by some
$h\in \pspdd$  then $g$ is also induced by $h$.
\end{lema} 

\begin{proof}  
Through any point $V\in \goke$ passes a $\P^{2(n-k)+2}$ of the first
kind. Indeed, it suffices to choose
$V^{k-1}\subset V$ of dimension $k-1$ and $H$ of
dimension $2n-k+2$ such that $V\subset H\subset (V^{k-1})^{\bot}$,
and to consider the linear space $L$ they define by \eqref{eq:first_kind}. 
Then $gL$ passes through $gV$ and the commuting square \eqref{eq:carre_g_tilde} 
shows then that for all  $V\in \goke$ and $V^{k-1}\in G(k-1,E)$ we have
$$
V^{k-1}\subset V  \donne \tilde g V^{k-1} = h( V^{k-1})  \subset gV.
$$  
Since the hyperplane  $ V^{k-1} \subset V$ is arbitrary, it follows that $gV=h(V)$.
\end{proof}
\end{numar}

\section{A Borel-Weil theorem for the odd symplectic group}

We set out to prove the Borel-Weil type theorem stated in the
introduction. We start by examining the preferred classes of modules
which appear as the spaces of global sections of line bundles on the
odd symplectic flag manifold $\fodd$.

\begin{numar}  \textbf{Shtepin's class of $\lietspdd$-modules.} \label{nr:shtepin}  
In \cite{Shtepin}, Shtepin deals with the problem of separating
multiple components when restricting simple $\liespddd$-modules to
$\liespd$. There is no semi-simple Lie algebra sitting between
$\liespddd$ and $\liespd$ which could be used as an intermediate step
in the reduction $\liespddd \downarrow \liespd$. Instead, Shtepin
considers the non-reductive intermediate Lie algebra $\lietspdd$ and
constructs in each simple $\liespddd$-module $V$ a filtration by
$\lietspdd$-modules 
$$
V_1 \subset \dots \subset V_p \subset V
$$
such that the factors $V_{i+1}/V_i$ are pairwise non-isomorphic and
multiplicity free $\liespd$-modules. 
It turns out that the $\lietspdd$-modules 
$V_{i+1}/V_i$ have nice properties, in particular as we will show
below,  they appear as
spaces of global sections of line bundles on $\fodd$. 
We now give a short account of  this construction.

We will use Shtepin's definition for the intermediate Lie algebra,
which is
\begin{equation}  \label{eq:shtepin_tspdd} 
\lietspdd = \{ X \in \liespddd \mid Xe_n=0 \}.
\end{equation}
This is not just a mere change of notation with respect to our
definition from \ref{nr:intermediate} since we
still want to
distinguish the Borel subalgebra $\liebddd$ of $\liespddd$ of upper
triangular matrices in the symplectic basis 
$ \{ e_0, \dots, e_{2n+1}  \} $ as giving the positive roots and thus
the notion of dominant and highest weights.
Of course, we can get back to our definition from \ref{nr:intermediate}
by conjugating by some element $g \in \spddd$ which sends $e_n$ to
$e_0$, but then we will also have to conjugate the Borel subalgebra
and consider as positive roots the roots of $g\liebddd g^{-1}$, which
is rather awkward. 
In accordance with the definition \eqref{eq:shtepin_tspdd} we now embed
$\liespd \subset \liespddd$ as the subalgebra
$$
\liespd=\{X\in\liespddd\mid Xe_n=0,\, Xe_{\bar{n}}=0 \}.
$$

Let  $\lietddd$ be the Cartan subalgebra of $\liespddd$ of diagonal
matrices in the 
symplectic basis $ \{ e_0, \dots, e_{2n+1}  \} $ and 
$\liebddd^-$ the Borel subalgebra opposite to $\liebddd$.
Denote $\lieumnsddd$, $\lieuplsddd$ the maximal nilpotent subalgebras
of $\liebddd^-$ and $\liebddd$.
Let $\lietd = \lietddd \cap \liespd$, 
 $\liebd=\liebddd\cap\liespd$,
respectively  $\lieu^{\pm}_{2n} =\lieu^{\pm}_{2n+2}  \cap\liespd$ be the
 corresponding subalgebras of $\liespd$.
Denote also
$\lieu^{\pm}_{2n+1}  =\lieu^{\pm}_{2n+2}  \cap\lietspdd$ and
$\liebdd=\liebddd\cap\lietspdd$.
Then $\liebdd= \lietd \oplus \lieuplsdd$, but note that $\liebdd$ is
not a Borel subalgebra of $\lietspdd$.
We have the decomposition 
\begin{equation} \label{eq:cartan_tspddd} 
\lietspdd =\lieumnsdd\oplus \lietd\oplus \lieuplsdd.
\end{equation}
Note that since we have
\begin{align*} 
\lieuplsdd&
  =\lieuplsd\oplus \bigoplus _{0\leq i \leq n-1} 
              \C X_{i\bar{n}} \oplus \C X_{n\bar n}  \\
\lieumnsdd&
  =\lieumnsd\oplus \bigoplus _{0\leq i \leq n-1} 
              \C X_{ni},
\end{align*}
the $\lietd$-weights of $\lieuplsdd$ are the positive roots of
$\liespd$ to which we add the $\e_i$, for $1 \leq i \leq n-1$ (and 0
corresponding to the element $X_{n\bar n}$) and the $\lietd$-weights
of $\lieumnsdd$ are their negatives (except for the 0 weight which
does not appear). Having this analogy between the decomposition
\eqref{eq:cartan_tspddd} and a Cartan decomposition is the reason to
choosing
the definition \eqref{eq:shtepin_tspdd} instead of our earlier one.

Let $V$ be a simple $\liespddd$-module. Shtepin constructs his
filtration in $V$ using an explicit model of $V$ as a space of
functions on the unipotent group $\uplsddd$ with Lie algebra
$\lieuplsddd$ due to {\v{Z}}elobenko (cf. \cite{Zhl}), which we won't
detail here. Following Shtepin, we call \emph{semi-maximal} 
the $\lietddd$-eigenvectors of $V$ which are killed by $\lieuplsd$.
These are exactly the highest 
weight 
vectors of the simple pieces of a
decomposition of $V$ as a $\liespd$-module.
Among the semi-maximal vectors we choose those which are killed by
the elements $X_{0n}, X_{1n}, \dots, X_{n-1\,n}$ of $\lieuplsddd$,  
which we call \emph{quasi-maximal}.
We order them
$v_1, \dots, v_p$ in a natural way coming from their explicit
expressions in the {\v{Z}}elobenko model and we define then the
filtration
$V_1 \subset \dots \subset V_p $ 
in
$V$ by
taking $V_i$ to be the $\lietspdd$-module generated by $v_1, \dots, v_i$.

The ordering $v_1, \dots, v_p$ is chosen such that for all $1 \leq i
\leq p-1$ we have
$\lieuplsdd v_{i+1} \in V_i$, so the image $\bar v_{i+1}$ of $v_{i+1}$
in $V_{i+1} / V_i$ is a vector killed by $\lieuplsdd$, ie a
\emph{maximal} vector with respect to $\lietspdd$.
The
factor $V_{i+1}/ V_i$ is a cyclic $\lietspdd$-module generated by the
maximal vector $\bar v_{i+1}$, and so, if $\mu $ is the
$\lietd$-weight of $\bar v_{i+1}$, $V_{i+1}/ V_i$ is a quotient of the
Verma-like module
\begin{equation}  \label{eq:Verma} 
U(\lietspdd) \tens _{U(\liebdd)}  \C _\mu 
\end{equation}
where $\C_\mu $  is the $\liebdd$-module of dimension 1 on which
$\liebdd$ acts through the weight $\mu $.
This already accounts for some properties of the $\lietspdd$-module
$V_{i+1}/V_i$
similar to properties of simple modules for semi-simple Lie algebras,
such as:
\begin{enumerate} 
\item $\mu $ is the highest weight of $V_{i+1}/V_i$, ie the weights of
  $V_{i+1}/V_i$ are of the form $\mu -\theta $ with $\theta $ in the
  semi-group generated by the positive roots of $\lietspdd$;
\item the Weyl group $W(\spd)$ acts on the weights of  $V_{i+1}/V_i$
  and for any $w \in \wspd$ the $w\mu $-eigenspace is of dimension 1;
\item any $\lietspdd$-module endomorphism of  $V_{i+1}/V_i$ is scalar,
  in particular  $V_{i+1}/V_i$ is indecomposable;
\item  $V_{i+1}/V_i$ has a lowest weight, which is $\bar w_0\mu $
  where $\bar w_0$ is the longest element of $\wspd$;
\item[etc.] 
\end{enumerate}
Actually, if $\lambda $ is the highest weight of the simple
$\liespddd$-module $V$, then the quasi-maximal vectors of $V$ (and so
the factors  $V_{i+1}/V_i$) are
parametrized by patterns $\mu  \to \lambda $, where $\mu $ is also an
$\liespddd$-dominant weight and where the notation $\mu  \to \lambda $
means 
$$
\lambda _1\geq\mu _1\geq\lambda _2
     \geq\cdots\geq 
           \mu _n\geq\lambda_{n+1}\geq \mu _{n+1} .
$$
The quasimaximal vector $q_{\mu ,\lambda }$ corresponding to the
pattern $\mu  \to \lambda $ is a $\lietddd$-eigenvector with weight
$$
\mu - (|\lambda| -  |\mu|) \e_n,
$$
in particular its $\lietd$-weight is $\res{\mu }{\lietd} $.
Let us denote by $V_\lambda $ the simple $\liespddd$-module with
highest weight $\lambda $. Shtepin proves then the following result:

\begin{thm} \label{thm:Shtepin} 
  The factor of the filtration of $V_\lambda $ generated by
  the image of the quasimaximal vector $q_{\mu , \lambda }$ is
  isomorphic to the $\lietspdd$-submodule of $V_\mu $ generated by the
  highest weight vector. 
\end{thm}

Denote by $V'_\mu $ this $\lietspdd$-module. Note that 
$V'_\mu $ does not depend only on 
$\res{\mu}{\lietd}$, 
as one might suspect given that it is a quotient of the
Verma module \eqref{eq:Verma}. In fact it is already so for
its structure as an $\liespd$-module, as Shtepin shows that the decomposition of
$V'_\mu $ into irreducible $\liespd$-modules is 
\begin{equation}  \label{eq:branching}  
\bigoplus_{\nu \to \mu } V''_{\nu},
\end{equation}
where $V''_\nu $ is the simple $\liespd$-module of highest weight 
$\nu$ and where the notation $\nu \to \mu$ now means
$$
\mu_1\geq\nu _1\geq\mu _2\geq\cdots\geq\mu _n\geq\nu _n\geq\mu _{n+1}.
$$

We will also need to consider the parabolic subalgebra
$$
\liep=\{X\in\liespddd\mid X(\C e_n)\subset\C e_n\}.
$$
We have
\begin{equation*} \label{eq:desc_liep}  
\liep     =\lietspdd  \oplus \C X_{nn}
\end{equation*}
and since $X_{nn} \in \lietddd$, the modules $V_i$ of the filtration
of a simple $\liespddd$-module $V$  are also $\liep$-modules.
Similarly, the  $\lietspdd$-submodule $V'_\mu $ of a simple
$\liespddd$-module $V_\mu$ generated by the highest weight vector is
a $\liep$-module. 
Then from the theorem \ref{thm:Shtepin} it follows that the factor of
the filtration of a simple $\liespddd$-module $V_\lambda $ generated by 
the image of the quasimaximal vector $q_{\mu , \lambda }$ is 
isomorphic to the $\liep$-module
$$
V'_\mu  \tens  \C_{-(|\lambda| -  |\mu|) \e_n}.
$$
\end{numar}

\begin{numar} \textbf{The $\liespdd$-modules of Proctor.} 
We continue to use the definition \eqref{eq:shtepin_tspdd} and
accordingly, by $\liespdd$ we mean the Lie algebra of the odd
symplectic group $\spdd$ associated to the odd symplectic form
$\res{\omega }{e_n^{\bot} }$. Recall from \ref{nr:tldo} that the restriction
morphism $g \mapsto \res{g}{e_n^{\bot} } $ induces a surjective
morphism $\liep \to \liespdd$ which restricts to a surjective morphism
$\lietspdd \to \liespdd \cap \liesldd$. The kernel of both these
morphisms is the one dimensional space $\C X_{n\bar n}$. 

In \cite{Proctor1}, Proctor uses an adapted form of the construction
of Weyl to define a special class of  $\liespdd$-modules. Let us recall this construction.

Let $\omega \in \Lambda ^2V^*$ be a generic skew-form on a complex
vector space $V$, and let $\lambda $ be a partition.
The $\GL(V)$-module of highest weight $\lambda $ coincides with the
Schur power
$S_\lambda V\subset V^{\tens d}$, where $d=|\lambda |$.
For $1\leq p<q\leq d$, denote $\varphi _{pq}$ the contraction by
$\omega $ on the indices $p$ et $q$
\begin{equation} \label{eq:phi_pq}  
\begin{array}{c} 
\varphi _{pq}: V^{\tens d}\vers V^{\tens d-2} \\
v_1\tens\cdots\tens v_n\longmapsto \omega (v_p,v_q)v_1\tens\cdots\tens \widehat{v_p}\tens\cdots\tens \widehat{v_q}\tens\cdots\tens v_n.
\end{array}
\end{equation}
Following \cite[\S17.3]{FultonHarris}, denote
$$
V^{\Spc{d}}=\bigcap_{p,q} \Ker \varphi_{pq}\subset V^{\tens d}
$$
the space of ``trace free''  $d$-tensors, and finally
$$
\ssl  V=S_{\lambda }V\cap V^{\Spc{d}}.
$$
When $V=\cd$,
if $\ell(\lambda )>n$ then  $\ssl  \cd=0$, and 
if $\ell(\lambda )\leq n$ then $\ssl  \cd$ is the simple $\liespd$-module
of highest weight $\lambda $.

In \cite{Proctor1}, Proctor studies the $\liespdd$-modules
$\ssl  \cdd$. He shows for example the following
(cf. \cite[theorem 2.1]{Proctor1}) :
\begin{prop} 
If $\ell(\lambda )>n+1$ then $\ssl  \cdd=0$. 
If $\ell(\lambda )\leq n+1$ then the  $\liespdd$-module  $\ssl  \cdd$ is indecomposable.
\end{prop}
He further describes weight bases for the modules $\ssl\cdd$ and gives
a formula for the $\lietd$-character of $\ssl\cdd$ similar to Weyl's formula.
As it turns out (cf. \cite[corollary 8.1]{Proctor1}), the decomposition
of $\ssl\cdd$ in simple $\liespd$-components is also given by the
formula  \eqref{eq:branching}. 
This is not a mere coincidence, since we have:

\begin{prop} \label{prop:shtepin_proctor} 
As a $\liep$-module via the morphism $\liep \to \liespdd$,  $\ssl\cdd$ 
is isomorphic to the $\liep$-module $V'_\lambda $.
\end{prop}

\begin{proof}
Let $\tE=\cddd$ and $\tldo \in \Lambda ^2\tE^*$ be the symplectic form
on $\tE$ for which the basis $ \{e_0, \dots, e_{2n+1}  \}$ is
symplectic. Denote $E= e_n^{\bot}$ and $\omega = \res{\tldo}{E}$.   
We identify the simple $\liespddd$-module $V_\lambda $ of highest
weight $\lambda $ with $\ssl \tE$. The highest weight vector
$v_\lambda \in V_\lambda $ coincides in this case with the image in
$S_\lambda \tE$ of the tensor
$$
e^{\lambda }= 
         e_0^{\tens \lambda _0}\tens e_1^{\tens \lambda _1}
                  \tens\cdots\tens 
         e_n^{\tens \lambda _n}    \in \tE^{\tens |\lambda |}.
$$
But since $e^{\lambda } \in E^{\tens |\lambda |}$, we have 
$v_\lambda \in S_\lambda  E \subset S_\lambda \tE$. 
On the other hand, $v_\lambda \in E^{\Spc{|\lambda |} } $  since
clearly the contractions \eqref{eq:phi_pq} defined by $\omega $ are
just the restrictions to $E^{\tens|\lambda |} $ of the contractions
\eqref{eq:phi_pq} defined by $\tldo$. 
Therefore $v_\lambda  \in \ssl E$.

Since $\liep$ acts on $\ssl E$ through the restriction morphism
$\liep \to \liespdd$, it follows that the $\liep$-module $V'_\lambda $
generated by $v_\lambda $ is actually contained in $\ssl E$. But
$V'_\lambda $ and $\ssl E$ have the same dimension, since they are
isomorphic as $\liespd$-modules. Therefore they must coincide.
\end{proof}
\end{numar}

\begin{numar} \textbf{Borel-Weil for the odd symplectic group.}  \label{sect:borel_weyl} 
Let $G$ be a semi-simple complex Lie group, $B \subset G$ a Borel subgroup
and $T \subset B$ a maximal torus. Denote $W$ the associated Weyl
group.
Let $e^1$ be the fixed point in the flag variety $G/B$ of the
opposite Borel subgroup $B^-$. Denote then, for $w \in W$, 
$e^w = w(e^1)$ and $X^w$ the opposite Schubert variety which is the closure of
the $B^-$-orbit of $e^w$. For a weight $\lambda $ let
$L_\lambda $ be the line bundle on $G/B$ whose fiber at the point
$e^1$ is the $B^-$-module $\C_\lambda $. By the Borel-Weil theorem, if
$\lambda $ is dominant then
$H^0(G/B, L_\lambda ) \simeq V_\lambda$, where $V_\lambda$ is the
simple $G$-module of highest weight $\lambda$. 
Denote the restriction of $L_\lambda$ to the Schubert variety $X^w$
also by $L_\lambda$.
Then (cf. \cite[\S 14.19]{Jantzen}) the restriction morphism
$$
H^0(G/B, L_\lambda ) \vers H^0(X^w, L_\lambda)
$$
is surjective and 
\begin{equation}  \label{eq:Demazure} 
H^0(X^w, L_\lambda )^* \simeq \Span(B^- \cdot v_{-w\lambda })
    \subset  V_\lambda ^*
\end{equation}
where $v_{-w\lambda } \in V_\lambda ^*$ is a vector of weight
$-w\lambda$ (unique up to a scalar factor) and 
$\Span(B^- \cdot v_{-w\lambda })$ is the $B^-$-submodule it generates.

Denote also $e_1$ the fixed point in $G/B$ of $B$, and, for $w \in W$,
$X_w$ the Schubert variety which is the closure of the $B$-orbit of
$e_w=w(e_1)$. Then $X_w = w_0 X^w$, where $w_0$ is  the longest
element of $W$. Conjugating by $w_0$ in \eqref{eq:Demazure}  
we get 
\begin{equation}  \label{eq:Xw}  
H^0(X_w, L_\lambda )^* \simeq \Span(B \cdot v_{-w_0w\lambda })
    \subset  V_\lambda ^*.
\end{equation} 
  
Take now $G= \spddd$, $\bddd \subset \spddd$ the Borel subgroup of upper
triangular matrices in the symplectic basis 
$ \{ e_0, \dots, e_{2n+1} \}$, 
$\tddd \subset \bddd$  the maximal torus of diagonal
matrices in the basis $ \{ e_0, \dots, e_{2n+1} \}$,
and denote $E=e_0^{\bot}$. 
Let
$w=\wz \in \wspddd$ so that, as in \ref{sect:cell_sch_fodd},
$\foe  = X_{\wz}$. 
Note that  the longest element $w_0$ of the Weyl group $\wspddd$ acts
as $-1$ and $V_\lambda ^* \simeq V_\lambda $.
From \eqref{eq:Xw} we get
\begin{equation}  \label{eq:BW_fodd} 
H^0(\foe, L_\lambda )^* \simeq \Span(\bddd \cdot v_{w\lambda })
    \subset  V_\lambda.
\end{equation}
Since $\foe$ is $P$-stable, where $P$ is the parabolic subgroup
preserving the line $\C e_0$, 
we have 
$\Span(\bddd \cdot v_{w\lambda }) = \Span(P \cdot v_{w\lambda })$
and
this isomorphism 
is also an isomorphism of $P$-modules.

Denote also $E' = e_n^{\bot}$ and consider $\fo(E')$ which embeds in
$\foddd$ as 
$$ 
\{ (V_1\subset\dots\subset V_{n+1}) \in\foddd 
           \mid V_{n+1}\subset E'\}.
$$
Let $P'$ be the parabolic subgroup preserving the line $\C e_n$ and
$\liep'$ its Lie algebra (which in \ref{nr:shtepin} has been denoted
$\liep$). 

\begin{prop}  
As a $\liep'$-module, $H^0( \fo(E'), L_\lambda )^*$ 
is isomorphic to Shtepin's $V'_\lambda$.
\end{prop}

\begin{proof} 
Since $w(n)=0$, we have $\fo(E') = w^{-1} \foe$ and $P'=w^{-1}Pw$, so
conjugating by $w$ in \eqref{eq:BW_fodd}, we get
$H^0( \fo(E'), L_\lambda )^* \simeq \Span(P' \cdot v_\lambda)=V'_\lambda$.
\end{proof} 

\begin{cor} \label{cor:BW_shtepin} 
As $\liep'$-modules,  $H^0( \fo(E'), L_\lambda )^* \simeq \ssl E'$.
\end{cor}
\begin{proof} This follows immediately from
  \ref{prop:shtepin_proctor}.
\end{proof}

Let now $T_i$ be the tautological bundle of rank $i$ on the symplectic
flag manifold $\foddd$, whose fiber at a point 
$(V_1\subset\dots\subset V_{n+1})$ is $V_i$. The restriction 
$\res{T_i}{\foe}$ is then the rank $i$ tautological bundle on 
$\foe$ and we denote it also by $T_i$.
We have $T_{i+1}/T_i = L_{-\e_i}$
for  $1\leq i\leq n$ 
and
$T_1=L_{-\e_0}$
so 
$$
L_\lambda  = 
{T_1^*}^{\tens \lambda _0}    \tens
    {(T_2/T_1)^*}^{\tens \lambda _1}    \tens\cdots\tens 
    {(T_{n+1}/T_n)^*}^{\tens \lambda _n}.
$$

\begin{thm} \label{thm:BW}  
Let $\lambda= (\lambda_0 \geq\cdots\geq  \lambda_n \geq0)$
  be a partition and $L_\lambda$ the line bundle on $\fo(E)$
$$
L_\lambda  = 
{T_1^*}^{\tens \lambda _0}    \tens
    {(T_2/T_1)^*}^{\tens \lambda _1}    \tens\cdots\tens 
    {(T_{n+1}/T_n)^*}^{\tens \lambda _n}.
$$
Then as $\spdd$-modules, 
$H^0(\foe,L_\lambda)^* \simeq \ssl E$.
\end{thm}

\begin{proof}
First denote $\Sp'_{2n+1}$ the odd symplectic group which is the image
of $P'$ via the restriction morphism $g \mapsto \res{g}{E'}$. Then the
modules in \ref{cor:BW_shtepin} are also $\Sp'_{2n+1}$-modules. 
We have $\spdd=w \Sp'_{2n+1} w^{-1}$ and the result follows from 
\ref{cor:BW_shtepin} once we
observe that $\ssl E = w \ssl E'$ as subspaces of $\ssl \cddd$.  
\end{proof} 

Note that the statement of theorem 
\ref{thm:BW} is independent of the embedding $\foe \subset \foddd$.
\end{numar}

\bibliography{bib} 
\bibliographystyle{amsalpha} 

\end{document}